%% file: Sto_maxcut_v3_arxiv.tex
\documentclass[11pt]{article}
\usepackage[margin=1in]{geometry}

\usepackage{multirow}
\usepackage[utf8]{inputenc}
\usepackage{amsmath}
\usepackage{amsfonts}
\usepackage{color}
\usepackage{latexsym}
\usepackage{tablefootnote}
\usepackage{algorithm}
\usepackage{algorithmic}
\usepackage{hyperref}

\newcommand*{\Var}{\mathsf{Var}}
\newcommand*{\QEDA}{\hfill\hbox{\vrule width1.0ex height1.0ex}}

\include{def}

\begin{document}

 \title{
 Multi-cut
 stochastic approximation methods 
\\
for solving stochastic	
	convex composite optimization
}
\date{May 21, 2025 (1st revision: March 1, 2026; 2nd revision: October 22, 2026)}
\author{
		Jiaming Liang \thanks{Goergen Institute for Data Science and Artificial Intelligence (GIDS-AI) and Department of Computer Science, University of Rochester, Rochester, NY 14620 (email: {\tt jiaming.liang@rochester.edu}). This work was partially supported by GIDS-AI seed funding and AFOSR grant FA9550-25-1-0182.}
        \qquad 
        Renato D.C. Monteiro \thanks{School of Industrial and Systems
		Engineering, Georgia Institute of
		Technology, Atlanta, GA, 30332-0205.
		(email: {\tt renato.monteiro@isye.gatech.edu} and {\tt hzhang906@gatech.edu}). This work
			was partially supported by AFOSR Grants FA9550-25-1-0131.}\qquad 
	Honghao Zhang \footnotemark[2]
	}
	
\maketitle
\begin{abstract}


This paper considers
the stochastic convex composite optimization  problem and
presents
multi-cut stochastic approximation (SA) methods for solving it, whose models in expectation overestimate its objective function.
The multi-cut model obtained by taking the maximum of a finite number of linearizations of the stochastic objective function provides a biased estimate of the objective function, with the error being uncontrollable.
Instead, our proposed SA method uses models obtained by taking the maximum of a finite number of one-cut models, i.e., suitable convex combinations of linearizations of the stochastic objective function.
It is shown that the proposed methods achieve nearly optimal convergence rate and have computational performance comparable,   and sometimes superior, to other SA-type methods.\\



		\par {\bf Keywords.} stochastic convex composite optimization, stochastic approximation, multi-cut method, optimal complexity bound.
		\\
		
		{\bf AMS subject classifications.} 
		49M37, 65K05, 68Q25, 90C25, 90C30, 90C60.
	\end{abstract}
\section{Introduction}

This paper considers
the
stochastic	
	convex composite optimization (SCCO) problem
\begin{equation}\label{eq:ProbIntro}
	\phi_{*}:=\min \left\{\phi(x):=f(x)+h(x): x \in \R^n\right\}
	\end{equation}
	where 
\begin{equation}\label{pbint2}
	f(x)=\mathbb{E}_{\xi}[F(x,\xi)].
	\end{equation}
	The following conditions are assumed:
	i) $f, h: \R^{n} \rightarrow \R\cup \{ +\infty \} $ are proper closed convex functions
	such that
	$ \dom h \subseteq \dom f $;
	ii) for almost every $\xi \in \Xi$,
a convex functional oracle $F(\cdot,\xi) :\dom h \to \R$ and
a stochastic subgradient
oracle $s(\cdot,\xi):\dom h \to \R^n$ satisfying
$s(x,\xi) \in \partial F(x,\xi)$
for every $x \in \dom h$ are available;
	and iii) for every $x \in \dom h$, $\Exp{\|s(x,\xi)\|^2} \le M^2$ for some $M\in \R_+$. The main goal of this paper is i) to present a multi-cut stochastic approximation (SA) method
 whose cutting-plane models are pointwise maximum of suitable one-cut models, and ii) to establish a global convergence rate $\tilde {\cal O}(1/\sqrt{j})$ at iteration $j$.

An important application of SCCO is two-stage stochastic programming. For example, consider
\[
\min_{x\in X}\;
c^\top x+\mathbb{E}_\xi[Q(x,\xi)], \quad \mbox{where } Q(x,\xi)= \min_{y\geq 0} \left\{ q(\xi)^\top y: T(\xi)x+W(\xi)y=r(\xi) \right\}.
\]
This is a special case of \eqref{eq:ProbIntro}-\eqref{pbint2} with $F(x,\xi)=c^\top x+Q(x,\xi)$ and $h=\delta_X$ being the indicator function of $X$.
Under standard recourse assumptions \cite{shadenrbook}, the value and a
subgradient of $Q(\cdot,\xi)$ can be obtained by solving the second-stage problem.

SA methods for solving \eqref{eq:ProbIntro} are iterative schemes that use a sequence of stochastic models to approximate the objective function $\phi$.   
Given an initial point $z_0$, many SA methods solve a sequence of prox subproblems given by
\begin{equation}\label{def:xj1} 
        z_{j} =\underset{u\in \R^n}\argmin
	    \left\lbrace  \Gamma_j^\lam(u):=
	    \Gamma_j(u) +\frac{1}{2\lam}\|u- z_0 \|^2 \right\rbrace
\end{equation}
where $\Gamma_j(\cdot)$ is some stochastic model for $\phi(\cdot)$. All methods developed in the literature for solving the SCCO problem \eqref{eq:ProbIntro} approximate the  objective function $\phi$ by stochastic models $\Gamma_j$ satisfying 
\begin{equation}\label{ineq:common}
\Exp{\Gamma_j(u)} \le  \phi(u), \quad \forall u \in \dom h,
\end{equation}
a key condition which plays a fundamental role in their complexity analysis.

The goal of this paper is to develop multi-cut SA methods whose models $\Gamma_j(\cdot)$ are the pointwise maxima of stochastic linear approximations of $\phi(\cdot)$. 
Stochastic models $\Gamma_j$ constructed in this manner will generally violate \eqref{ineq:common}.
An obvious example of a  multi-cut model $\Gamma_j$ that fulfills the above description is
\begin{equation}\label{multicut}
\Gamma_j(\cdot) = \max_{1\le i\le j}\{\ell(\cdot,z_{i-1};\xi_{i-1})\},
\end{equation}
where 
$\ell(\cdot,x;\xi):= h(\cdot) + [ F(x,\xi)+\inner{s(x,\xi)}{\cdot-x}]$ is the (composite) linearization of $\phi$ at the point~$x$. 
(Here, ``multi-cut'' refers to the pointwise maximum of multiple stochastic cutting-plane models of the expected objective, rather than the scenario-wise cuts used in multi-cut Benders decomposition.)
However, this model may significantly violate \eqref{ineq:common} since all that is known is that $\Exp{\Gamma_j(\cdot)}= \phi(\cdot) + {\cal O}(\sqrt{j})$ for every $j$ (see \eqref{ineq:badbd} below), which makes it an unsuitable ingredient for the development of convergent multi-cut SA methods.

{\bf Literature Review.}
The first SA method for solving \eqref{eq:ProbIntro} was introduced by Robbins and Monro \cite{monro51}  under the assumption that $h \equiv 0$ and, for almost every $\xi \in \Xi$, 
$F(\cdot,\xi)$ is smooth and convex. Since then, a variety of SA methods have been developed, including stochastic (sub)gradient methods \cite{nemjudlannem09,nemyud78,polyak90,polyakjud92}, stochastic mirror descent \cite{dang2015stochastic,nemyud78}, stochastic accelerated gradient methods \cite{ghadimi2012optimal,ghadimi2013optimal,lan2012optimal}, stochastic dual averaging (DA) methods of \cite{nesterov2009primal,xiao2009dual}, and the SCPB 
 method of \cite{liang2024single}.
All of these SA methods use single-cut models
$\Gamma_j(\cdot)$ of the form
\begin{equation}\label{onecut}
     \Gamma_j(\cdot) := \sum_{i=1}^j \alpha_i^j\,\ell(\cdot,z_{i-1};\xi_{i-1}),
\end{equation}
 for some $ \alpha_i^j \in [0,1] $, $1\le i \le j$, satisfying $\sum_{i=1}^j \alpha_i^j=1$.
 Consequently, a common characteristic of these SA methods is that all
$\Gamma_j$'s are \emph{single-cut} approximations, and hence satisfy \eqref{ineq:common}.
Except for stochastic DA and SCPB methods, the aforementioned SA methods employ the single-cut model in \eqref{onecut} with $\alpha^j_j=1$, meaning that only the most recent stochastic (sub)gradient is retained and historical information is discarded. 
In contrast, stochastic DA \cite{nesterov2009primal,xiao2009dual} and SCPB \cite{liang2024single} adopt the single-cut model in \eqref{onecut} where the coefficients $\alpha_1^j,\ldots,\alpha_j^j$ are all positive.
As a result, these methods construct approximation models $\Gamma_j$ using all the history up to the current iteration.

Instead of directly computing $\Gamma_j$ using \eqref{onecut},
 both SCPB and stochastic DA update the model $\Gamma_j$ recursively via the formula $\Gamma_j(\cdot) = \beta_{j}\Gamma_{j-1}(\cdot)+(1-\beta_j)\ell(\cdot,x_{j-1};\xi_{j-1})$ for some $\beta_j \in (0,1)$. This update requires only the previous model $\Gamma_{j-1}$ and the most recent sample $\xi_{j-1}$, avoiding the need to store all past samples explicitly. 
In contrast, constructing the $j$-th multi-cut model $\Gamma_j$ in \eqref{multicut} requires access to all past random samples $(\xi_0, \ldots, \xi_{j-1})$, and hence entails significantly greater storage than the aforementioned SA methods.
As the iteration count $j$ grows, this storage demand can become substantial, particularly in real-world applications such as two-stage stochastic programming, where each realization $\xi$ often includes a large data matrix.

 Methods that construct cutting-plane models as in \eqref{multicut}, referred to as \emph{multi-cut} algorithms, have also been extensively studied in the stochastic programming literature, 
such as stochastic decomposition (SD) and sample average approximation (SAA) methods.
SD methods were first proposed in \cite{higle1991stochastic}.
Building on Ruszczyński's regularized deterministic decomposition algorithm \cite{rus1986ruglarized},
a sampling-based regularized SD variant was subsequently developed and further extended in \cite{higle1994finite,higle1996stochastic}. Most of the papers on SD methods (if not all) focus on two-stage or multi-stage stochastic programs \cite{gangammanavar2021stochastic,higle1991stochastic,higle1994finite,higle1996stochastic,liu2020asymptotic,sen1994network,sen2014multistage} and,
to the best of our knowledge, the authors are not aware of any paper that considers SD methods in the context of the SCCO problem \eqref{eq:ProbIntro}.  More recently, \cite{liu2020asymptotic} established  convergence rate results for SD methods applied to two-stage stochastic quadratic programs. SD methods allow the random vector $\xi$ to have an arbitrary distribution, either discrete or continuous, and iteratively construct a sequence of piecewise linear approximations expressed as the pointwise maximum of affine functions.
The computation of the $j$-th approximation model $\Gamma_j$
requires that
$j$ (or at most $n+3$) random samples be available, and
hence may
demand a substantial amount of storage as $j$ becomes large.
Although SD methods fall within the class of multi-cut algorithms, these approximations all satisfy the condition \eqref{ineq:common}.

Another prominent multi-cut approach is the SAA method \cite{kleywegt2002sample,shapiro1991asymptotic,shapiro2003monte,verweij2003sample}, which approximates the expected objective \eqref{pbint2} by an empirical average $\sum_{i=0}^{J-1}F(\cdot, \xi_i)/J$ for a large i.i.d. sample $(\xi_0, \ldots,\xi_{J-1})$ of $\xi$ at the beginning of the method and then applies certain deterministic methods to minimize the empirical average.
The SAA method must retain all random samples $(\xi_0, \ldots,\xi_{J-1})$ in memory in order to generate new subgradients of the empirical average.
Again, this can result in substantial
storage requirements when the sample size $J$ is large. 
Among SAA variants, the L-shaped method \cite{vwets} and the regularized L-shaped method \cite{rus1986ruglarized} are two well-known methods that employ multi-cut models, due to their practical performance in solving large-scale two-stage
stochastic programming problems. 
From a nonsmooth optimization perspective, these methods can be interpreted as cutting-plane and proximal bundle methods \cite{Kelley1960,lemarechal1975extension,lemarechal1978nonsmooth,mifflin1982modification,wolfe1975method} applied to the empirical average $\sum_{i=0}^{J-1}F(\cdot, \xi_i)/J$.  

 {\bf Contributions:} 
The main contribution of this paper is the development of new multi-cut SA methods for solving the  SCCO problem \eqref{eq:ProbIntro} with nearly optimal convergence rate guarantees.

Motivated by the classical multi-cut model \eqref{multicut} and the single-cut model \eqref{onecut}, we first design a single-phase multi-cut SA method, termed the single-phase max-one-cut (S-Max1C) method, whose models  $\Gamma_j$  are constructed as maximum of one-cut SA models of the form \eqref{onecut} initialized at different iterations. Unlike previously studied SA methods, the SA models $\Gamma_j$ generated by S-Max1C do not necessarily satisfy condition \eqref{ineq:common}, but 
instead the relaxed condition
\[
\mathbb{E}[\Gamma_j(u)] - \phi(u) = \mathcal{\tilde O}\!\left(\frac{1}{\sqrt{j}}\right), \quad \forall u \in \dom \phi,
\]
which is sufficient to derive a nearly optimal convergence rate for S-Max1C. We further develop a multi-phase variant via a warm-start strategy and prove that it preserves the same nearly optimal convergence rate. Our analysis accommodates very general sampling assumptions: the random sample $\xi$ may follow an arbitrary distribution, either discrete or continuous. Finally, from a computational perspective, the max model $\Gamma_j$ constructed by S-Max1C requires only $\mathcal{O}(\log j)$  random samples of $\xi$, which is substantially smaller than the ${\cal O}(j)$ samples required to build the $j$-th model in the SD methods. Thus, S-Max1C demands substantially less storage than SD methods, while incurring only a modest memory overhead compared to the single-cut SA methods discussed in the second paragraph of the Literature Review.

{\bf Organization of the paper.} Subsection~\ref{subsec:DefNot}  presents basic definitions and notation used throughout the paper.
Section~\ref{sec:framework} formally describes the assumptions on the SCCO problem \eqref{eq:ProbIntro} and
presents a stochastic cutting plane (S-CP) framework which is used to analyze some important instances
contained in it.
Section~\ref{sec:method} presents the S-Max1C method of the S-CP framework and establishes its convergence rate bound.
Section~\ref{sec:multistage} provides a multi-phase version of the S-Max1C method and its convergence analysis.
Section~\ref{sec:proof} presents the deferred proof of the main result of Section~\ref{sec:framework}.
Section~\ref{sec:numerics} presents computational results to illustrate the efficiency of our proposed methods.
Section~\ref{sec:conclusion} presents some concluding remarks and possible extensions.
Finally, Appendix~\ref{app:technical} contains technical results used in our analysis.
   
    \subsection{Basic definitions and notation}\label{subsec:DefNot}
    
    Let $\mathbb{N}_{++}$ denote the set of positive integers.
    The sets of real numbers, non-negative, and positive real numbers are denoted by $\R$, $\R_+ $, and $\R_{++}$, respectively. 
	Let $\R^n$ denote the standard $n$-dimensional Euclidean 
	space equipped with  inner product and norm denoted by $\left\langle \cdot,\cdot\right\rangle $
	and $\|\cdot\|$, respectively. Throughout the paper, $\log(\cdot)$ denotes the natural logarithm.

	Let $\psi: \R^n\rightarrow (-\infty,+\infty]$ be given. The effective domain of $\psi$ is denoted by
	$\dom \psi:=\{x \in \R^n: \psi (x) <\infty\}$ and $\psi$ is proper if $\dom \psi \ne \emptyset$.
 For $\varepsilon \ge 0$, the \emph{$\varepsilon$-subdifferential} of $ \psi $ at $z \in \dom \psi$ is defined as
    \[
    \partial_\varepsilon \psi (z):=\left\{ s \in\R^n: \psi(u)\geq \psi(z)+\left\langle s,u-z\right\rangle -\varepsilon, \forall u\in\R^n\right\}.
    \]
	The subdifferential of $\psi$ at $z \in \dom \psi$, denoted by $\partial \psi (z)$, is by definition the set  $\partial_0 \psi(z)$.
	Moreover, for some scalar $\mu \ge 0$, a proper function $\psi: \R^n\rightarrow (-\infty,+\infty]$ is said to be $\mu$-convex if
	$$
	\psi(\alpha z+(1-\alpha) u)\leq \alpha \psi(z)+(1-\alpha)\psi(u) - \frac{\alpha(1-\alpha) \mu}{2}\|z-u\|^2
	$$
	for every $z, u \in \dom \psi$ and $\alpha \in [0,1]$. Denote the set of all proper lower semicontinuous convex functions by $\bConv{n}$.

\section{A Stochastic Cutting Plane Framework}\label{sec:framework}
This section contains three subsections. The first one describes the assumptions made on problem \eqref{eq:ProbIntro}. The second one presents the motivation for our work. 
The third one presents an S-CP framework, which is used in Section \ref{sec:method} to analyze some specific
cases contained in it. 

\subsection{Assumptions on the SCCO problem}\label{Assumptions}
	
Let $\Xi$ denote the
	support of random 
	vector $\xi$ and assume that the following conditions on the SCCO problem \eqref{eq:ProbIntro} hold:
	\begin{itemize}
\item[(A1)]
$f$ and $ h$ are proper closed convex functions satisfying
		$\dom f \supseteq \dom h$;
		
\item[(A2)] for almost every $\xi \in \Xi$,
a convex functional oracle $F(\cdot,\xi) :\dom h \to \R$ and
a stochastic subgradient
oracle $s(\cdot,\xi):\dom h \to \R^n$ satisfying
\begin{equation}\label{eq:unbias}
    f(x) = \Exp{F(x,\xi)}, \quad s(x,\xi) \in \partial F(x,\xi) \quad \forall x \in \dom h,
\end{equation}
 are available;
 
		\item[(A3)]
		$M := \sup \{ \Exp{\|s(x,\xi)\|^2}^{1/2} : x \in \dom h \} < \infty$;
 	\item[(A4)]
		the set of optimal solutions $X _*$ of
	 \eqref{eq:ProbIntro} is nonempty.
	\end{itemize}

We need some definitions and basic facts in the analysis of the paper. For every $x \in \dom h$, let
\begin{equation}\label{def:f'}
f'(x) := \Exp{s(x,\xi)} \in \partial f(x),
\end{equation}
where the inclusion directly follows from \eqref{eq:unbias}.
Moreover, for every
$\xi \in \Xi$ and $x \in \dom h$, let
\begin{equation}\label{def:Phi}
	\Phi(\cdot,\xi)=F(\cdot,\xi)+h(\cdot), \quad \ell(\cdot,x;\xi) := F(x,\xi)+ h(\cdot)+ \IInner{s(x,\xi)}{\cdot-x}.
	\end{equation} 
Condition (A3) and Jensen's inequality imply that for every $x \in \dom h$,
\begin{equation}\label{ineq:fp}
 \|f'(x)\| = \|\Exp{s(x,\xi)}\| \le \Exp{\|s(x,\xi)\|} \le \left(\Exp{\|s(x,\xi)\|^2}\right)^{1/2} \le M.
\end{equation}
Also, the definitions of $\ell$ and $\Phi$ in \eqref{def:Phi} and
the fact that $F(\cdot,\xi)$ is convex due to (A2)
imply that for every $u \in \dom h$,
\begin{equation}\label{ineq:linear}
   \ell(u,x;\xi) \le \Phi(u,\xi). 
\end{equation}
Hence, $\ell(\cdot,x;\xi)$ is
a stochastic composite linear approximation of
$\phi(\cdot)$
in the sense that its expectation is a true composite linear approximation of $\phi(\cdot)$.

\subsection{Motivation for this work} \label{subsec:motivation}
This subsection outlines how S-Max1C differs from  other state-of-the-art methods developed  in the literature for solving \eqref{eq:ProbIntro}, or special cases of it.

Like S-Max1C, the comparison only considers methods, namely,  SD~\cite{higle1991stochastic}, RSA~\cite{nemjudlannem09}, SCPB \cite{liang2024single}, and DA~\cite{nesterov2009primal,xiao2009dual}, which solve a finite sequence  of prox 
subproblems 
(possibly only one)
\begin{equation}\label{def:turesub}
\underset{u\in \R^n}\min
	    \left\lbrace  \phi(u)+\frac{1}{2\lam}\|u- z_0 \|^2 \right\rbrace.
\end{equation}
Each subproblem is uniquely determined by a pair $(\lam,z_0)$, which varies from one subproblem to another. For a fixed pair $(\lam,z_0)$, all these methods  solve \eqref{def:turesub} by solving a sequence of subproblems as in \eqref{def:xj1} where the $\Gamma_j$'s are models constructed using some update rules. In the following, we discuss how the models $\Gamma_j$ generated by these methods compare to the ones used by S-Max1C.

 We first introduce some definitions. For some pair of indices $(k,j)$ such that $1 \le k \le j$, consider 
 the collection $ \mathcal{C}_{k,j}$ of functions given by
\begin{equation}\label{eq:defPhiij}
\sum_{i=k}^j \alpha_i\,\ell(\cdot,z_{i-1};\xi_{i-1}),
\end{equation}
where $\ell(\cdot,\cdot;\cdot)$ is as in \eqref{def:Phi} and the scalars $\alpha_i \ge 0$ for every $i=k,\ldots,j$ 
satisfy $\sum_{i=k}^j \alpha_i=1$. If $\Psi(\cdot) \in\mathcal{C}_{k,j}$, then 
\(
\mathbb{E}[\Psi(\cdot)] \le \phi(\cdot)
\) 
  due to relation \eqref{ineq:linear} and the fact that \(\mathbb{E}[\Phi(\cdot,\xi_{i-1})]=\phi(\cdot)\) for every $i=k,\ldots,j$. 

RSA can be viewed as a method for approximately solving \eqref{def:turesub} in the following sense: it performs a single iteration  
that consists of solving \eqref{def:xj1}
with $\lam$ sufficiently small and $\Gamma_1$  set to $\ell(\cdot,z_{0};\xi_{0})$, and hence with $\Gamma_1 \in \mathcal{C}_{1,1}$.
On the other hand, DA and SCPB approximately solve \eqref{def:turesub} by performing multiple iterations,
each of which solves a subproblem as in  \eqref{def:xj1} for a suitably generated model $\Gamma_j \in \mathcal{C}_{1,j}$. Hence, all the models mentioned above satisfy the 
 condition
\eqref{ineq:common} in view of the last remark in the previous paragraph.

On the other hand,
the $j$-th model \(\Gamma^{\mathrm{our}}_j\) constructed by our method S-Max1C has the form
\begin{equation}\label{def:Gammaour}
\Gamma^{\mathrm{our}}_j
=\max\bigl\{L_k^j: k\in B_j \bigr\},
\end{equation}
for some index set  $B_j$  such that
$|B_j| = {\cal O}(\log j)$ and functions $ L_k^j(\cdot) \in {\cal C}_{k,j}$ for every $k \in B_j$.
Despite the fact that 
$\mathbb{E}[L_{k}^j(\cdot)] \le \phi(\cdot) $ for every $k \le j$,
the latter inequality does not ensure that $\Gamma_j=\Gamma_j^{our}$ satisfies the 
 condition
\eqref{ineq:common}.
One of the main contributions  of this paper is to show that $\Gamma_j^{our}$ satisfies the relaxed inequality
\[
\mathbb{E}\bigl[\Gamma^{\mathrm{our}}_j(\cdot)\bigr]
\le \phi(\cdot)+\tilde {\cal O}\left(\frac{1}{\sqrt{j}}\right)
\]
 as long as
$B_j$ is not too large, i.e.,
$|B_j| = {\cal O}(\log j)$ (see Proposition~\ref{prop:noise} below).

Moreover, it is shown in the remarks following the algorithm in Subsection 2.2 of \cite{higle1991stochastic} that, at iteration \(j\), SD constructs a max-model $\Gamma^{\mathrm{SD}}_j$ such that
$\Gamma^{\mathrm{SD}}_j $ is the pointwise maximum of $j$ affine functions
 $A_1,\ldots,A_j$
 all underneath the function $\Phi_1^j(\cdot):= \sum_{i=1}^j \,\Phi(\cdot,\xi_{i-1})/j$, i.e., $A_i(\cdot) \le \Phi_1^j(\cdot)$ for every $i=1,\ldots,j$ and
\begin{equation}\label{def:GammaSD}
\Gamma^{\mathrm{SD}}_j(\cdot)
=\max\bigl\{A_i(\cdot) : i=1,\ldots,j\bigr\}.
\end{equation} 
Hence,
\[
\Gamma^{\mathrm{SD}}_j(\cdot) \le \Phi_1^j(\cdot),
\qquad\text{and}\qquad
\mathbb{E}\bigl[\Gamma_j^{\mathrm{SD}}(\cdot)\bigr]
\le \mathbb{E}\bigl[\Phi_1^j(\cdot)\bigr]=\phi(\cdot),
\]
where the identity is due to the first relation in \eqref{eq:unbias}. 
In conclusion, among the five methods discussed above, our approach is the only one that allows its models to violate the  condition~\eqref{ineq:common}.
 Moreover, the fact that S-Max1C violates this condition benefits its computational performance as demonstrated by the computational results in Section~\ref{sec:numerics}.
 
Finally, while the SD model $\Gamma^{\mathrm{SD}}_j$ is the maximum of  ${\cal O}(j)$ functions,  the  S-Max1C model $\Gamma_j$ is the  maximum of only
${\cal O}(\log j)$  functions.

\subsection{S-CP framework}\label{subsec:SPB}

This subsection describes the S-CP framework
and states a result that provides a bound on the optimality gap of its iterates. 

We make some remarks about S-CP.
First, S-CP is a single-phase method, i.e., it solves a single prox subproblem as in \eqref{def:turesub}, and hence keeps the prox center the same throughout.
It differs from standard proximal point-type methods that solve a sequence of prox subproblems as in \eqref{def:turesub}, and hence are multi-phase methods.
 Second, S-CP is not an implementable algorithm but is rather a conceptual framework since the  flexible condition \eqref{def:generalGamma}
 does not specify $\Gamma_j(\cdot)$. 
Third, S-CP outputs a triple $(z_I,z_I^a,u_I)$ including the last iterate $z_I$, the average of iterates $z_I^a$, and the average of observed stochastic function values $u_I$.

Throughout our analysis, we let
\[
\xi_{[j]}=(\xi_0,\xi_1,\ldots,\xi_j) \quad
\forall
j=0,\dots ,I.
\]
Since steps 1-3 of S-CP imply that $z_{j}$
depends on $\xi_{[j-1]}$
but not on $\xi_{j}$, and
the latter is chosen to be independent of $\xi_{[j-1]}$, we conclude that 
$z_{j}$ is independent of
$\xi_{j}$, and hence that \eqref{eq:unbias} and \eqref{def:f'} imply that
\begin{equation}\label{eq:unbiasedzj}
f(z_j) = \Exp{ F(z_j,\xi_j)}, \quad
f'(z_j) = \Exp{ s(z_j,\xi_j)} \quad \forall j=1,\ldots,I.
\end{equation}

\begin{algorithm}[H]
\caption{S-CP framework}
\label{alg:SCP}

\begin{algorithmic}

\STATE \textbf{Input:} Scalar $\lam>0$, integer $I\ge 1$, and point
$z_0\in\dom h$.

\STATE \textbf{0.}
Generate a sample $\xi_0$ of $\xi$, set $j=1$, and
\begin{equation}
\label{def:beta}
\beta
=
\frac{I+1-\log(I+1)}
     {I+1+\log(I+1)};
\end{equation}

\STATE \textbf{1.}
Evaluate the stochastic functional and subgradient oracles of (A2)
at $(z_{j-1},\xi_{j-1})$ to obtain the pair
$\big(F(z_{j-1},\xi_{j-1}),s(z_{j-1},\xi_{j-1})\big)$
and the linearization
$\ell(\cdot,z_{j-1};\xi_{j-1})$ as in~\eqref{def:Phi};

\STATE \textbf{2.}
Choose $\Gamma_j\in\bConv{n}$ such that
\begin{align}
\label{def:generalGamma}
\Gamma_j(\cdot)
\ge
\begin{cases}
\ell(\cdot,z_{j-1};\xi_{j-1}),
    & \text{if } j=1,\\
\beta\Gamma_{j-1}(\cdot)
 +(1-\beta)\ell(\cdot,z_{j-1};\xi_{j-1}),
    & \text{otherwise};
\end{cases}
\end{align}

\STATE \textbf{3.}
Compute $z_j$ as in~\eqref{def:xj1} and
\begin{align}
\label{def:wj1}
z_j^a
=
\begin{cases}
z_j,
    & \text{if } j=1,\\
(1-\beta)z_j+\beta z_{j-1}^a,
    & \text{otherwise};
\end{cases}
\end{align}
Generate a sample $\xi_j$ of $\xi$, independently of
$\xi_0,\ldots,\xi_{j-1}$, and compute
\begin{align}
\label{def:u}
u_j
=
\begin{cases}
\Phi(z_1,\xi_1),
    & \text{if } j=1,\\
(1-\beta)\Phi(z_j,\xi_j)+\beta u_{j-1},
    & \text{otherwise};
\end{cases}
\end{align}

\STATE \textbf{4.}
If $j=I$, then \textbf{stop}; otherwise, set
$j\leftarrow j+1$ and go to Step~1.

\STATE \textbf{Output:} $(z_I,z_I^a,u_I)$.

\end{algorithmic}
\end{algorithm}

The following technical result reveals a simple relationship between $u_j$ and $\phi(z_j^a)$ in expectation.

\begin{proposition}
  For every $j \ge 1$, it holds that  
\begin{equation}\label{ineq:basic1}
    \Exp{u_j} \ge \Exp{\phi(z_j^a)} \ge \phi_* .
\end{equation}
\end{proposition}
\begin{proof}
The second inequality in \eqref{ineq:basic1} obviously follows from the definition of $\phi_*$ in \eqref{eq:ProbIntro} and the fact that $z_j^a \in \dom h$ for every $j \ge 1$. We next prove the first inequality in \eqref{ineq:basic1}  by induction on $j$. If $j=1$, then
it follows from
\eqref{def:wj1}, \eqref{def:u}, and \eqref{eq:unbiasedzj} that
\[
	\Exp{u_1} \overset{\eqref{def:u}}{=}
\Exp{\Phi(z_{1},\xi_{1})}
 \overset{\eqref{eq:unbiasedzj}}{=}
\Exp{\phi(z_{1})} 
 \overset{\eqref{def:wj1}}{=} \Exp{\phi(z_1^a)}.
 \]
Assume now that the first inequality in \eqref{ineq:basic1}  holds
for $j-1$ for some $j\ge 2$. This assumption, relations
\eqref{def:wj1}, \eqref{def:u}, and \eqref{eq:unbiasedzj}, 
and the fact that  $\phi$ is convex, then imply that
\begin{align*}
    \Exp{u_{j}}  &\overset{\eqref{def:u},\eqref{eq:unbiasedzj}} = (1-\beta) \Exp{\phi(z_{j})} + \beta \Exp{u_{j-1}}\overset{\eqref{ineq:basic1}}{\ge}  (1-\beta) \Exp{\phi(z_{j})} + \beta \Exp{\phi(z_{j-1}^a)}\\
     &\ge \Exp{\phi((1-\beta)z_{j} + \beta z_{j-1}^a)}  \overset{\eqref{def:wj1}}{=}\Exp{\phi(z_{j}^a)},
\end{align*}
and hence that
\eqref{ineq:basic1}  holds.
\end{proof}

We now present two concrete examples of the model $\Gamma_j$ satisfying \eqref{def:generalGamma}. These two models set $\Gamma_1(\cdot):= \ell(\cdot,z_0;\xi_0)$ and construct $\Gamma_j$ for $j \ge 2$
recursively as follows:
 \begin{itemize}
	    \item [(S1)] The {\bf one-cut} update scheme sets    
        $\Gamma_j(\cdot) := \beta \Gamma_{j-1}(\cdot) + (1-\beta)\ell(\cdot,z_{j-1};\xi_{j-1})$.
  It is easy to see that the models generated by this scheme are given by  
\begin{equation}\label{eq:Gammaj}
    \Gamma_j(\cdot) := \beta^{j-1}\ell(\cdot,z_{0};\xi_{0}) + (1-\beta) \sum_{i=2}^{j} \beta^{j-i}\ell(\cdot,z_{i-1};\xi_{i-1}).
\end{equation}
Moreover,  (A2), inequality \eqref{ineq:linear}, and the above identity,
imply that
\begin{equation}\label{ineq:NjS1}
    \Exp{ \Gamma_j(\cdot) } \stackrel{\eqref{ineq:linear}}\le  \Exp{\beta^{j-1} \Phi(\cdot,\xi_{0}) + (1-\beta) \sum_{i=2}^{j} \beta^{j-i} \Phi(\cdot,\xi_{i-1})}  = \phi(\cdot).
\end{equation}       
	   	    
	    \item [(S2)] The \textbf{multi-cut} update scheme  sets  
	        $ \Gamma_j(u) = \max \{ \Gamma_{j-1}(u),\ell(u,z_{j-1};\xi_{j-1})\}$. It is easy to see that  the models generated by this scheme are given by  
\begin{equation}\label{eq:multicut}
\Gamma_j(\cdot) = \max_{1\le i\le j}\{\ell(\cdot,z_{i-1};\xi_{i-1})\}.
\end{equation}
Moreover, it is shown in Lemma \ref{lem:Nj} that for every $u \in \dom h$,
\begin{equation}\label{ineq:badbd}
\Exp{ \Gamma_j(u)} - \phi(u) \le 2\sigma(u) \sqrt{j-1},
\end{equation}
where 
$\sigma(u)$ is defined as
\begin{equation}
    \label{def:sigma_u}
    \sigma(u):= \sqrt{\Exp{(\Phi(u,\xi)-\phi(u))^2}}.
\end{equation}
	\end{itemize}
For every $j \in \{1,\ldots,I\}$,  define the $j$-th model noise at $u \in \dom h$ as
\begin{equation}\label{def:Nj}
{\cal N}(u;\Gamma_j) = 
 \max \left\{ 0, \Exp{ \Gamma_j(u)} - \phi(u)\right\}.
\end{equation}
Clearly, ${\cal N}(u;\Gamma_j) \ge 0$ for any $u \in \dom h$.
 It follows from  \eqref{ineq:NjS1} that ${\cal N}(u;\Gamma_j)$ is always $0$ if $\Gamma_j$ is generated according to (S1) but,
 in view of \eqref{ineq:badbd},
 it might  grow as ${\Theta}(\sqrt{j})$   if $\Gamma_j$ is generated according to (S2).

We next state a result that provides a preliminary bound on $\Exp{ u_I}$. Its proof is postponed to Section~\ref{sec:proof} since it follows along similar (but simpler) arguments as those used in the analysis of \cite{liang2024single}.
For the following two results, we need the definitions below
 \begin{equation}\label{def:d0}
    x_* := \argmin_{x \in X_*} \|z_0-x\|, \quad d_0 := \|z_0-x_*\|.
\end{equation}

\begin{proposition}\label{lem:t1-2}
S-CP with input $(\lam,I,z_0)$ satisfies
\begin{equation}\label{keyrecu}
   \Exp{ u_I} - \phi_* \le \frac{44\lam M^2 \log (I+1)}{I+1} + \frac1{2\lam}\Exp{\|z_0 -x_*\|^2}- \frac{1 }{2 \lam}\Exp{\| z_I - x_*\|^2}+ {\cal N}(x_*;\Gamma_I). 
 \end{equation}

\end{proposition}

By suitably choosing the prox stepsize $\lam$, one obtains the following important specialization of Proposition~\ref{lem:t1-2}.

\begin{theorem}
\label{coro1}
Let  $D \ge d_0$ be an over-estimate of $d_0$ and $\lam$ be given by
\begin{equation}\label{eq:lamtau1} 
      \lambda= \frac{  D \sqrt{I+1}}{2 M \sqrt{ 22\log (I+1) }}.
\end{equation}
Then, S-CP with input $(\lam,I,z_0)$ satisfies
\begin{equation}\label{ineq:CS-Max1Conecut}
  \Exp{\phi( z_I^a)} -  \phi_* \le \Exp{u_I} -  \phi_*  \le \frac{2 M D\sqrt{ 22\log (I+1)}}{\sqrt{I+1}} +{\cal N}(x_*;\Gamma_I).     
\end{equation}
\end{theorem}

\begin{proof}
The first inequality in \eqref{ineq:CS-Max1Conecut} directly follows from \eqref{ineq:basic1}. Inequality  \eqref{keyrecu} and the definition of $d_0$ in \eqref{def:d0} imply that 
\[
    \Exp{u_I} - \phi_* \le  \frac{44\lam M^2\log (I+1)}{I+1} + \frac{d_0^2}{2\lam} + {\cal N}(x_*;\Gamma_I).
\]
Hence, the second inequality \eqref{ineq:CS-Max1Conecut} follows from the above inequality, the choice of $\lam$ in \eqref{eq:lamtau1} and the fact that $D \ge d_0$. 
\end{proof}

  It can be seen from \eqref{ineq:CS-Max1Conecut} that the first term on its right-hand side converges to $0$ at the rate of $ {\cal O}(1/\sqrt{I})$, and hence that the convergence rate of S-CP for one-cut bundle model is $ {\cal O}(1/\sqrt{I})$  since ${\cal N}(x_*;\Gamma_I) = 0$ in view of \eqref{ineq:NjS1} and \eqref{def:Nj}.

The next section considers an instance of S-CP that uses a bundle model $\Gamma_j$ lying between the two extreme cases (S1) and (S2) mentioned above.
Specifically, it chooses $\Gamma_j$ to be the maximum of one-cut models and shows that the noise
${\cal N}(u;\Gamma_j) = \tilde {\cal O}(1/\sqrt{j})$ for every $u \in \dom h$.

\section{The Max-one-cut Method}\label{sec:method}

As illustrated in \eqref{ineq:badbd}, the multi-cut scheme \eqref{eq:multicut} lacks control over the noise in function value, i.e., ${\cal N}(\cdot;\Gamma_I)$.
In this section, we consider a special instance of S-CP
where $\Gamma_j$ is obtained
by taking maximum of one-cut (Max1C) models.

For every $j\ge k \ge 1$, the one-cut model that starts at iteration $k$ and ends at iteration $j$ is defined as
\begin{equation}\label{def:L}
    L_k^j(\cdot) := \sum_{i=k}^{j} \beta_i^j\ell(\cdot,z_{i-1};\xi_{i-1}),
\end{equation}
where
\begin{equation}\label{def:betaij}
    \beta_k^j = \beta^{j-k}, \quad \beta_i^j = (1-\beta)\beta^{j-i},  \quad\forall  i \in \{k+1,\ldots, j\},
\end{equation}
and $\beta$ is as in \eqref{def:beta}.

Specifically, let $\{1\} \subseteq B \subseteq \{1,\ldots, I\}$
denote the set of indices at which the computation of a one-cut model is started
and, for a given $j \in \{1,\ldots, I\}$, let
\begin{equation}\label{def:Bj}
    B_j = \{ k \in B : k \le j \}
\end{equation}
denote the initial iteration indices of the one-cut models constructed up to the $j$-th iteration. 
The $j$-th bundle model is then defined as
\begin{equation}\label{def:Gammaj}
    \Gamma_j (\cdot) = \max_{ k \in B_j} L_k^j(\cdot).
\end{equation}
We note that this bundle model includes the one-cut scheme (S1) and the multi-cut scheme (S2) as special cases.
If $B=\{1\}$, then for every $1\le j\le I$, $B_j=\{1\}$ and $\Gamma_j (\cdot) = L_1^j(\cdot)$, i.e., $\Gamma_j$ is the same as \eqref{eq:Gammaj}, and hence it reduces to the one-cut scheme (S1).
If $B=\{1,\ldots, I\}$ and $\beta=1$, then for every $1\le j\le I$, $B_j=\{1, \ldots, j\}$ and $L_k^j (\cdot) =\ell(\cdot,z_{k-1};\xi_{k-1})$ in view of \eqref{def:L}, so $\Gamma_j$ is the same as \eqref{eq:multicut}, and hence it reduces to the multi-cut scheme (S2).

We state below the instance of the S-CP framework,
referred to as S-Max1C, where $\Gamma_j(\cdot)$ is selected as in
\eqref{def:Gammaj}.
However, instead of constructing
$\Gamma_j(\cdot)$ directly from its definition in \eqref{def:Gammaj}, step 2 uses a recursive formula for building $\Gamma_j$ from $\Gamma_{j-1}$ and the most recent linearization
$\ell(\cdot,z_{j-1};\xi_{j-1})$.
This recursive formula, which is shown in Lemma \ref{ineq:recurL} below, immediately implies that
$\Gamma_j$ in \eqref{def:Gammaj} satisfies condition \eqref{def:generalGamma} imposed by the S-CP framework.\\





\begin{algorithm}[htbp]
\caption{S-Max1C method}
\label{alg:SMax1C}

\begin{algorithmic}

\STATE \textbf{Input:} Scalar $\lam>0$, integer $I\ge 2$, set $B$ such that
$\{1\} \subseteq B \subseteq \{1,\ldots,\lfloor I/2\rfloor\}$,
and initial point $z_0\in\dom h$.

\STATE \textbf{0.}
Same as Step~0 of Algorithm~\ref{alg:SCP};

\STATE \textbf{1.}
Same as Step~1 of Algorithm~\ref{alg:SCP};

\STATE \textbf{2.}
Compute
\begin{equation}
\label{eq:keyobs}
\Gamma_j(\cdot)
=
\begin{cases}
\ell(\cdot,z_0;\xi_0),
    & \text{if } j=1,\\[1mm]
(1-\beta)\ell(\cdot,z_{j-1};\xi_{j-1})
+\beta\max\left\{
\Gamma_{j-1}(\cdot),
\ell(\cdot,z_{j-1};\xi_{j-1})
\right\},
    & \text{if } j\in B\setminus\{1\},\\[1mm]
(1-\beta)\ell(\cdot,z_{j-1};\xi_{j-1})
+\beta\Gamma_{j-1}(\cdot),
    & \text{otherwise};
\end{cases}
\end{equation}

\STATE \textbf{3.}
Same as Step~3 of Algorithm~\ref{alg:SCP};

\STATE \textbf{4.}
If $j=I$, then \textbf{stop}; otherwise, set
$j\leftarrow j+1$ and go to Step~1.

\STATE \textbf{Output:} $(z_I,z_I^a,u_I)$.

\end{algorithmic}
\end{algorithm}

Now we make some remarks about S-Max1C. First, the first one-cut model used within the max model $\Gamma_j(\cdot)$ starts at the first iteration. The condition that
$B \subseteq \{1,\ldots, \lfloor I/2 \rfloor\}$ ensures that
the max-one-cut models are constructed with
one-cut models that start during its first half iterations and thus $B_j = B$ for all $j > I/2$ where $B_j$ is as in \eqref{def:Bj}. Second, two practical choices of $B$ are: i) $B=\{1\}$; and   ii) $B$ is equal to the set of all the powers of 2 less than or equal to $I/2$, namely, 
\begin{equation}\label{B:2_i}
    B = \left\{2^i: 0 \le 2^i \le \frac{I}2\right\}.
\end{equation}
Now we state the main convergence result for S-Max1C.
\begin{theorem}\label{thm:main1}
Assume that $\sigma(x_*) < \infty$.
Let a positive
integer $ I \ge 2$ and a set $B$ such that 
$\{1\} \subseteq B \subseteq \{1,\ldots, \lfloor I/2 \rfloor\}$ be given, and $\lam$ be as in \eqref{eq:lamtau1}.
Then, S-Max1C with input $(\lam,I,B,z_0)$ satisfies
\begin{equation}\label{ineq:CS-Max1C}
    \Exp{\phi( z_I^a)} - \phi_* \le \Exp{u_I} - \phi_* \le \frac{2\sqrt{\log (I+1)}}{\sqrt{I+1}}
    \left[\sqrt{22} M D+2 \sigma(x_*)\sqrt{|B|-1}\right].
\end{equation}
\end{theorem}
Now we make some remarks about Theorem \ref{thm:main1}. First, when $B = \{1\}$, $\Gamma_j(\cdot)$ in \eqref{eq:keyobs} reduces to the one-cut model (S1), inequality \eqref{ineq:CS-Max1C} then implies that the convergence rate of S-CP for one-cut bundle model is $ \tilde {\cal O}(1/\sqrt{I})$. Second, inequality \eqref{ineq:CS-Max1C} shows that the convergence rate for S-Max1C is $\tilde {\cal O}(1/\sqrt{I})$ if $|B| = {\cal O}(\log I)$.

The remainder of this section is devoted to the proof of Theorem \ref{thm:main1}. 
The analysis begins with the following result showing that S-Max1C is an instance of the S-CP framework.

\begin{lemma}\label{ineq:recurL}
   The bundle model $\Gamma_j(\cdot)$ defined in \eqref{def:Gammaj}  satisfies the recursive formula in step 2 of S-Max1C. Moreover, S-Max1C is an instance of the S-CP framework.
\end{lemma}

\begin{proof}
When $j=1$, it follows from \eqref{eq:keyobs} that   $\Gamma_1 =\ell(\cdot,z_0;\xi_0)$, which is the same as \eqref{def:Gammaj} with $j=1$. Throughout the remaining proof, we assume $j \ge 2$.
First, the definitions of $L_k^j$ and $\beta_i^j$ in \eqref{def:L} and \eqref{def:betaij}, respectively, imply that for every $j\ge k \ge 1$,
    \begin{equation} \label{eq:Lkj}
	L_k^j(\cdot) =  \left\{\begin{array}{ll}
		   (1-\beta)\ell(\cdot,z_{j-1};\xi_{j-1}) + \beta L_k^{j-1}(\cdot),  & \text { if } j \ge k+1,
     \\ 
		   \ell(\cdot,z_{k-1};\xi_{k-1}),  & \text { if } j=k.
		    \end{array}\right. 
    \end{equation}
It is also easy to see from the definition of $B_j$ in \eqref{def:Bj} that
\[
B_j =  \left\{\begin{array}{ll}
		   B_{j-1} \cup \{j\},  & \text { if } j \in B,
     \\ 
		    B_{j-1} ,  & \text { if } j\notin B.
		    \end{array}\right. 
\]
    If $j \notin B$, then it follows from the definition of $B_j$ in \eqref{def:Bj} that $j\ge k+1$ for every $k \in B_j$.
    Using this observation, relation \eqref{eq:Lkj}, the definition of $\Gamma_j$ in \eqref{def:Gammaj}, and the fact that $B_j = B_{j-1}$, we have
\begin{align}
    \Gamma_{j}(\cdot) &\stackrel{\eqref{def:Gammaj}}= \max_{k \in B_j }L_k^j(\cdot) \stackrel{\eqref{eq:Lkj}}= (1-\beta)\ell(\cdot,z_{j-1};\xi_{j-1}) + \beta \max_{ k \in B_{j-1}} L_k^{j-1}(\cdot) \nn \\
&\stackrel{\eqref{def:Gammaj}}= (1-\beta)\ell(\cdot,z_{j-1};\xi_{j-1}) + \beta \Gamma_{j-1}(\cdot). \label{eq:G1}
\end{align}
If $j \in B$, then $B_j = B_{j-1} \cup \{j\}$. It thus follows from the definition of $\Gamma_j$ in \eqref{def:Gammaj} and relation \eqref{eq:Lkj} that
\begin{equation}\label{eq:max}
    \Gamma_{j}(\cdot) \stackrel{\eqref{def:Gammaj}}= \max_{k \in B_j }L_k^j(\cdot) = \max_{k \in B_{j-1} \cup \{j\} }L_k^j(\cdot) = \max\left\{\max_{ k \in B_{j-1}} L_k^j(\cdot), L_j^j(\cdot) \right\}.
\end{equation}
Note that $j\ge k+1$ for every $k \in B_{j-1}$ in view of the definition of $B_j$ in \eqref{def:Bj}.
Using this observation and relation \eqref{eq:Lkj}, we have
\begin{align*}
    \max_{ k \in B_{j-1}} L_k^j(\cdot) &\stackrel{\eqref{eq:Lkj}}= (1-\beta)\ell(\cdot,z_{j-1};\xi_{j-1}) + \beta \max_{ k \in B_{j-1}} L_k^{j-1}(\cdot) \\
    &\stackrel{\eqref{def:Gammaj}}=(1-\beta)\ell(\cdot,z_{j-1};\xi_{j-1}) + \beta \Gamma_{j-1}(\cdot).
\end{align*}
Plugging the above equation and the formula for $L_j^j(\cdot)$ in \eqref{eq:Lkj} into \eqref{eq:max}, we obtain
\[
\Gamma_{j}(\cdot) = \max\left\{(1-\beta)\ell(\cdot,z_{j-1};\xi_{j-1}) + \beta \Gamma_{j-1}(\cdot),\ell(\cdot,z_{j-1};\xi_{j-1}) \right\}.
\]
Therefore, \eqref{eq:keyobs} holds due to the above identity and \eqref{eq:G1}.

Finally, we prove that S-Max1C is an instance of the S-CP framework. It suffices to show that $\Gamma_j$ in \eqref{eq:keyobs} satisfies step 2 of S-CP.
Clearly, the recursive formula of  $\Gamma_j$ in \eqref{eq:keyobs} implies that $\Gamma_j \in \bConv{n} $ and for every $u \in \dom h$,
\[
\Gamma_j(u) \ge \beta \Gamma_{j-1}(u)  + (1-\beta)\ell(u,z_{j-1};\xi_{j-1}).
\]
Therefore, $\Gamma_j$ in \eqref{eq:keyobs} satisfies \eqref{def:generalGamma}, and the second claim of the lemma is proved.
\end{proof}

Recall that Lemma \ref{lem:Nj} shows that the $j$-th model noise ${\cal N}(\cdot;\Gamma_j)$ in  \eqref{def:Nj} for $\Gamma_j$ generated according to the multi-cut scheme in \eqref{eq:multicut} is ${\cal O}(\sqrt{j})$. Hence,
 the last (i.e., the $I$-th) model noise satisfies ${\cal N} (\cdot;\Gamma_I)= {\cal O}(\sqrt{I})$.

 Our goal in the next two results is to show that the last model noise for S-Max1C
 is much smaller than the one above, i.e.,
 it satisfies
 ${\cal N} (\cdot;\Gamma_I)= \tilde {\cal O}(1/\sqrt{I})$
  as long as $|B|={\cal O}(\log I)$.
 Before establishing this fact in Proposition~\ref{prop:noise} below, we first discuss some properties about the convex combination of the stochastic functions $\Phi(\cdot,\xi_{i-1})$,
$i=k,\ldots,j$, given by
 \begin{equation} \label{def:Q}
	Q_k^j(\cdot) = \sum_{i=k}^j \beta_i^j \Phi(\cdot,\xi_{i-1})
    \end{equation}	
    with $\beta_i^j$ as in \eqref{def:betaij}.
    First,  for every $j\ge k \ge 1$ and $u \in \dom h$, we have
\begin{equation}\label{eq:exp1}
        L_k^j(u) \le Q_k^j(u), \quad \Exp{Q_k^j(u)} = \phi(u).
    \end{equation}
  Indeed, the inequality in \eqref{eq:exp1} immediately follows from relation \eqref{ineq:linear} and the definitions of $L_k^j$ and $Q_k^j$ in \eqref{def:L} and \eqref{def:Q}, respectively.
The identity in \eqref{eq:exp1} follows from the definition of $Q_k^j$ in \eqref{def:Q} and the fact that $\Exp{\Phi(\cdot,\xi)}=\phi(\cdot)$.
 Second, the following technical result shows that the variance of $Q_k^I$ is $\tilde {\cal O}(1/I)$ as long as $k$ is not too close to $I$.

\begin{lemma}\label{lem:basicLQ}
For every $I\ge 2$, $k \le \lfloor I/2 \rfloor$, and $u \in \dom h$, we have 
 \begin{equation}\label{ineq:QkI}
     \Exp{\left (Q_{k}^I(u)-\phi(u) \right)^2 } \le \frac{4\log (I+1)}{I+1}\sigma^2(u)
 \end{equation}
 where $Q_k^I$ is as in \eqref{def:Q}.
\end{lemma}

\begin{proof}
Fix $u \in \dom h$.
Define $Z_i = \beta_i^I\left(\Phi(u,\xi_{i-1}) - \phi(u)\right)$ for $i =k,\ldots,I$.  
Using the definitions of $Q_k^I(\cdot)$ and $\beta_i^j$   in \eqref{def:Q} and \eqref{def:betaij}, respectively, the fact that $\{\xi_{i-1}\}_{k\le i\le I}$ are independent, and Lemma~\ref{lem:ind}, we have
 \begin{align*}
     &\Exp{\left(Q_{k}^I(u)-\phi(u)\right)^2}
\overset{\eqref{def:Q}}= \Exp{\left(\sum_{i=k}^{I}Z_i\right)^2} = \sum_{i=k}^I \Exp{Z_i^2}\\
& = \left(\beta_k^I\right)^2\Exp{  \left(\Phi(u,\xi_{k-1})-\phi(u)\right)^2} + \sum_{i=k+1}^{I} \left(\beta_i^I\right)^{2} \Exp{\left(\Phi(u,\xi_{i-1})-\phi(u)\right)^2}\\
&\overset{\eqref{def:betaij}}=\beta ^{2(I-k)}\Exp{  \left(\Phi(u,\xi_{k-1})-\phi(u)\right)^2} + (1-\beta)^2 \sum_{i=k+1}^{I} \beta^{2(I-i)} \Exp{\left(\Phi(u,\xi_{i-1})-\phi(u)\right)^2}.
\end{align*} 
The above relation and the definition of $\sigma(\cdot)$ in \eqref{def:sigma_u} imply that
\begin{equation}\label{ineq:beta}
   \Exp{\left(Q_{k}^I(u)-\phi(u)\right)^2} \overset{\eqref{def:sigma_u}}= \left(\beta^{2(I-k)} + (1-\beta)^2 \sum_{i=k+1}^{I} \beta^{2(I-i)}\right) \sigma^2(u) 
    \le \left(\beta^I + \frac{1-\beta}{1+\beta}\right) \sigma^2(u),
\end{equation}
where the inequality is due to the facts that $k \le \lfloor I/2 \rfloor$ and $\sum_{i=k+1}^{I} \beta^{2(I-i)} \le 1/(1-\beta^2)$.
Using \eqref{def:beta} and Lemma \ref{lem:tauI} with $C=I+1$, we have
\begin{equation}\label{obs3}
\frac{1-\beta}{1+\beta} \stackrel{\eqref{def:beta}}= \frac{\log (I+1)}{I+1}, \quad \beta^{I} \stackrel{\eqref{ineq:beta-1}}\le 3 \beta^{I+1} \stackrel{\eqref{ineq:beta-1}}\le \frac{3}{I+1} \le \frac{3\log (I+1)}{I+1},
\end{equation}
where the last inequality is due to $\log (I+1) \ge 1$ for every $I \ge 2$.
Finally, we conclude that \eqref{ineq:QkI} immediately follows from \eqref{ineq:beta} and \eqref{obs3}.
\end{proof}

The following proposition establishes the key technical result ${\cal N} (\cdot;\Gamma_I)= \tilde {\cal O}(1/\sqrt{I})$. Its proof relies on a technical result in the appendix, namely Lemma \ref{lem:variance}.

\begin{proposition}\label{prop:noise}
For every $u \in \dom h$, we have
\begin{equation}\label{ineq:GammaI}
{\cal N}(u;\Gamma_I) \le 4\sigma(u) \sqrt{|B|-1}\ \frac{\sqrt{\log (I+1)}}{\sqrt{I+1}}, 
\end{equation}
where $\sigma(u) $ is as in \eqref{def:sigma_u} and ${\cal N}(u;\Gamma_I)$ is as in \eqref{def:Nj}.
\end{proposition}
\begin{proof}
  Fix $u \in \dom h$. Relation \eqref{eq:exp1} with $j=I$ and inequality \eqref{ineq:QkI} imply that the random variables $\{X_k\}_{k \in B}$ and $\{Y_k\}_{k \in B}$, and the scalar $\sigma_X$, defined as
\[
X_k = Q_k^I(u) - \phi(u), \quad Y_k=L_k^I(u)-\phi(u),
\quad  \sigma_X = \frac{2\sqrt{\log (I+1)}}{\sqrt{I+1}}\sigma(u),
\]
satisfy \eqref{eq:condition}.
Hence, it follows from the conclusion of 
Lemma \ref{lem:variance} and
the definition of $\Gamma_I$ in \eqref{def:Gammaj} that
\[
\Exp{\Gamma_I(u)-\phi(u)} \overset{ \eqref{def:Gammaj}}= \Exp{ \max_{k \in B} \left\{  L_k^I(u)-\phi(u) \right \} } \stackrel{\eqref{ineq:E-max}}\le 4\sigma(u) \sqrt{|B|-1}\
\frac{\sqrt{\log (I+1)}}{\sqrt{I+1}} . 
\]
Finally, inequality \eqref{ineq:GammaI} follows from the above inequality and the definition of ${\cal N}(u;\Gamma_I)$ in \eqref{def:Nj}. 
\end{proof}

\vgap

We are now ready to prove Theorem \ref{thm:main1}.

\vgap

\noindent
	{\bf Proof of Theorem \ref{thm:main1}} Since S-Max1C is an instance of the S-CP framework (see Lemma~\ref{ineq:recurL}), Theorem~\ref{coro1} holds for S-Max1C.
Therefore, Theorem \ref{thm:main1} immediately follows from Theorem~\ref{coro1} and Proposition~\ref{prop:noise} with $u=x_*$.
\QEDA

\section{The Multi-phase Max1C Method}\label{sec:multistage}
This section presents a multi-phase version of S-Max1C, denoted as M-Max1C,  with a warm-start approach. Each phase consists of one call to S-Max1C, and the output of a phase is used as the initial point of the subsequent phase.

We describe the multi-phase version below in Algorithm~\ref{alg:MMax1C}.





\begin{algorithm}[H]
\caption{M-Max1C method}
\label{alg:MMax1C}

\begin{algorithmic}

\STATE \textbf{Input:} Scalar $\lam>0$, integers $I\ge 2$ and $L\ge 1$, set $B$ such that
$\{1\}\subseteq B\subseteq\{1,\ldots,\lfloor I/2\rfloor\}$,
and initial point $z_0\in\dom h$.

\STATE \textbf{0.}
Set $l=1$ and $x_0=z_0$;

\STATE \textbf{1.}
Compute
\[
(x_l,y_l,w_l)
=
\text{S-Max1C}(\lam,I,B,x_{l-1});
\]

\STATE \textbf{2.}
If $l<L$, then set $l\leftarrow l+1$ and go to Step~1;
otherwise, compute
\begin{equation}
\label{eq:output}
y_L^a
=
\frac{1}{L}\sum_{l=1}^L y_l,
\qquad
w_L^a
=
\frac{1}{L}\sum_{l=1}^L w_l,
\end{equation}
and \textbf{stop}.

\STATE \textbf{Output:} $(y_L^a,w_L^a)$.

\end{algorithmic}
\end{algorithm}

We now make some observations about M-Max1C.
First, the index $l$ counts the number of phases, and hence the final $l$, namely $L$,  is the total number of phases of M-Max1C. Clearly, M-Max1C with $L=1$ reduces to S-Max1C.
Second, every phase of M-Max1C calls S-Max1C in step~1 with input $x_{l-1}$ set to be the output of
the previous phase. Third, every call to S-Max1C in step~1 performs $I$ iterations, and hence the total number of iterations performed by M-Max1C is $LI$.

The following result establishes a convergence rate bound on the optimality gap for the final average iterate $y_L^a$ obtained in  \eqref{eq:output}, which holds 
for any choice of stepsize $\lam$.

\begin{proposition}
         \label{cor:cmplx1-practical1}
Assume that $\sigma(x_*) < \infty$. The output $(y^a_L,w^a_L)$ of 
M-Max1C
satisfies \begin{equation}\label{eq:general2}
      \Exp{\phi(y_L^a)} - \phi_* \le \Exp{w_L^a} - \phi_* \le  \frac{44\lam M^2 \log (I+1)}{I+1}  + \frac{d_0^2}{2\lam L} + 4\sigma(x_*) \sqrt{|B|-1}\ \frac{\sqrt{\log (I+1)}}{\sqrt{I+1}},
\end{equation}
where $(\lam,I,L,B)$ is as described in its  input, $d_0$ and $x_*$ are as in \eqref{def:d0}.
\end{proposition}

\begin{proof}
Throughout the proof, we keep $x_*$ fixed as defined in \eqref{def:d0} with the original initial point $z_0$. Note that the bound \eqref{keyrecu} holds for any fixed $x_*\in X_*$, since its proof does not use the closest-point property in \eqref{def:d0}.
Since the $l$-th phase of M-Max1C calls 
S-Max1C in its step1, it follows from Proposition~\ref{lem:t1-2} with 
$(z_0,z_I,u_I) = (x_{l-1},x_l,w_l)$ that
\begin{equation}\label{eq:keyrec}
    \Exp{w_l} - \phi_* 
\overset{\eqref{keyrecu}}\le  \frac{44\lam M^2 \log (I+1)}{I+1} +\frac{1}{2\lam}\Exp{\| x_{l-1} - x_*\|^2  } - \frac{1}{2\lam}\Exp{\| x_l - x_*\|^2  } + {\cal N}(x_*;\Gamma_I^l).
\end{equation}
 Summing \eqref{eq:keyrec} from $l=1$ to $l=L$,  dividing the resulting inequality by $L$, and using the definition of $d_0$ in \eqref{def:d0}, we have 
 \[
\frac{1}{L}\sum_{l=1}^L \Exp{w_l} - \phi_* \le  \frac{44\lam M^2 \log (I+1)}{I+1}
 +\frac{d_0^2}{2\lam L} + \frac{1}{L}\sum_{l=1}^L{\cal N}(x_*;\Gamma_I^l).
 \]
The above inequality, the definition of $w_L^a$ in \eqref{eq:output}, and Proposition \ref{prop:noise} with $u=x_*$ then imply that the second inequality in \eqref{eq:general2} holds. Finally, we prove the first inequality in \eqref{eq:general2}. Using the definition of $w_L^a$ in \eqref{eq:output} and the inequality \eqref{ineq:basic1}, we have
\[
\Exp{w_L^a} \overset{\eqref{eq:output}} = \frac{1}{L}\sum_{l=1}^L\Exp{w_l}\overset{\eqref{ineq:basic1}}\ge \frac{1}{L}\sum_{l=1}^L\Exp{\phi(y_l)} \overset{\eqref{eq:output}}\ge \Exp{\phi(y_L^a)},
\]
where the last inequality is due to the convexity of $\phi$ and the definition of $y_L^a$ in \eqref{eq:output}. Hence, the first inequality in \eqref{eq:general2} also holds.
\end{proof}

By properly choosing the stepsize $\lam$, one obtains the following specialization of Proposition~\ref{cor:cmplx1-practical1}.

\begin{theorem}
         \label{cor:cmplx1-practical2}
Assume that $\sigma(x_*) < \infty$. Let positive
integers $L$ and $I$ be given and  $D \ge d_0$ be an over-estimate of $d_0$,
define
\begin{equation}\label{eq:theta1-practical2}
    \lambda= \frac{  D \sqrt{I+1}}{2 M \sqrt{22 L \log (I+1) }}.
\end{equation}
Then, 
M-Max1C with input $(\lam,I,L,B)$
satisfies:
\begin{equation}\label{ineq:Dh2}
    \Exp{\phi( y_L^a)} - \phi_* \le \Exp{w_L^a} - \phi_* \le  \frac{2\sqrt{\log (I+1)}}{\sqrt{I+1}}
    \left[ \frac{\sqrt{22} M D }{\sqrt{L}}+ 2\sigma(x_*)\sqrt{|B|-1} \right] . 
\end{equation}
\end{theorem}

\begin{proof}
 This statement follows from Proposition~\ref{cor:cmplx1-practical1} with $\lam$ in \eqref{eq:theta1-practical2} and the fact that $D \ge d_0$.
\end{proof}

 Now we make some remarks about Theorem \ref{cor:cmplx1-practical2}. If $B=\{1\}$, \eqref{ineq:Dh2} implies that the convergence rate for M-Max1C is $ \tilde {\cal O}(1/\sqrt{LI})$, which matches the optimal  $ {\cal O}(1/\sqrt{LI})$ convergence rate up to a logarithmic factor.
More generally for $|B| >1$, \eqref{ineq:Dh2} implies that the convergence rate of M-Max1C  is  $ \tilde {\cal O}(1/\sqrt{LI})$
whenever
\[
L = {\cal O}\left(\frac{M^2 D^2}{\sigma(x_*)^2(|B|-1)}\right).
\]
Second, M-Max1C with $B=\{1\}$ is closely related to the SCPB
method of \cite{liang2024single}. However, in contrast to M-Max1C, which can arbitrarily choose
a constant length $I$ for its phases, SCPB computes its cycle (or phase) lengths using two  rules that
yield variable cycle lengths.

\section{Proof of Proposition~\ref{lem:t1-2}} \label{sec:proof}
This section presents the proof for Proposition~\ref{lem:t1-2}.
The first lemma presents some useful relationships needed for this section.

\begin{lemma} \label{lemfirstserious} For every
$j \in \{ 1,\ldots,I\}$, we have 
\begin{equation}\label{eq:basic}
    \Exp{\phi(z_j) -\ell (z_j,z_{j-1};\xi_{j-1})} \le  2M \sqrt{\Exp{\|z_j-z_{j-1}\|^2}}. 
\end{equation}
\end{lemma}

\begin{proof}
It follows from the definition of $\ell(\cdot,x;\xi)$ in 
\eqref{def:Phi} and the first identity in \eqref{eq:unbiasedzj} that
\begin{align*}
    \Exp{\phi(z_j) -\ell (z_j,z_{j-1};\xi_{j-1})} & = \Exp{f(z_{j}) - f(z_{j-1}) - \inner{s(z_{j-1},\xi_{j-1})}{z_{j}-z_{j-1}}} \\
    &\stackrel{\eqref{def:f'}}\le \Exp{\inner{f'(z_j) - s(z_{j-1},\xi_{j-1})}{z_j-z_{j-1}}},
\end{align*}
where the inequality follows from \eqref{def:f'}.
Applying the Cauchy-Schwarz inequality for random vectors (i.e., $\Exp{\inner{X}{Y}} \le \sqrt{\Exp{\|X\|^2}} \sqrt{\Exp{\|Y\|^2}}$), we obtain
\begin{equation}\label{ineq:CS}
    \Exp{\phi(z_j) -\ell (z_j,z_{j-1};\xi_{j-1})}
    \le \sqrt{\Exp{\|f'(z_j) - s(z_{j-1},\xi_{j-1})\|^2}} \sqrt{\Exp{\|z_j-z_{j-1}\|^2}}.
\end{equation}
Using the triangle inequality, the fact that $(a+b)^2\le 2(a^2+b^2)$, Assumption (A3), and \eqref{ineq:fp}, we have
\[
    \Exp{\|f'(z_j) - s(z_{j-1},\xi_{j-1})\|^2}
    \le 2\Exp{\|f'(z_j)\|^2 + \|s(z_{j-1},\xi_{j-1})\|^2} \stackrel{\eqref{ineq:fp}}\le 4M^2.
\]
Inequality \eqref{eq:basic} then follows from plugging the above inequality into \eqref{ineq:CS}.
\end{proof}

 The next technical result introduces a key quantity, namely, scalar $\alpha_j$ below,  and
provides a useful recursive 
relation for it. 

\begin{lemma}\label{lem:102}
For every
$j \in \{ 1,\ldots,I\}$, define
\begin{equation}
    \alpha_j:= \Exp{u_j - \Gamma_j^\lam(z_j)} - \frac{2\lam M^2(1-\beta)}{\beta}
\label{def:delta}
\end{equation}
where $\lambda>0$ is the prox stepsize input to the S-CP framework.
Then, the following statements hold:
\begin{itemize}
    \item [a)] we have
    \begin{equation} \label{ineq:alpha1}  
    \alpha_{1}\le  2 \lam M^2;
    \end{equation}
    \item[b)] for every $ 2 \le j \le I$, we have
        \begin{equation}\label{ineq:alphaj}
            \alpha_{j} \le \beta \alpha_{j-1}.
        \end{equation}
\end{itemize}
\end{lemma}

\begin{proof}
a) 
Relations \eqref{def:delta},   \eqref{def:u}, \eqref{def:generalGamma}, \eqref{eq:unbiasedzj}, and \eqref{eq:basic}, all with $j=1$,  imply that
\begin{align*}
        \alpha_{1} & \stackrel{\eqref{def:delta}}{\le} \Exp{ u_1 - \Gamma_1^{\lam}(z_{1}) }     \stackrel{\eqref{def:u}}{=} \Exp{\Phi(z_1,\xi_1) - \Gamma_1(z_{1}) -\frac{1}{2\lam}\|z_1-z_0\|^2  }   \\       &\stackrel{\eqref{def:generalGamma},\eqref{eq:unbiasedzj}}\le \Exp{\phi(z_1) -\ell(z_1,z_0;\xi_0) -\frac{1}{2\lam}\|z_1-z_0\|^2}  \\
        & \overset{\eqref{eq:basic}}\le 2M\sqrt{\Exp{ \left\|z_1-z_0\right\|^2}}-\frac{1}{2\lam}\Exp{\left\|z_1-z_0 \right\|^2} \le 2\lam M^2,
\end{align*}
where the  last inequality is due to the fact that $-a^2+ 2ab\le  b^2 $ with $a = \sqrt{\Exp{\|z_1-z_0\|^2}}/\sqrt{2\lam}$ and $b =\sqrt{2\lam} M$.

b) Let $2 \le j \le I$ be given.
It follows from the definition of  $\Gamma_j^\lam$ in  \eqref{def:xj1}, relation \eqref{def:generalGamma}, and the fact that $\beta <1$  that
\begin{align}   
\Gamma_j^\lam(z_j)-(1-\beta) &\ell(z_j,z_{j-1};\xi_{j-1}) \stackrel{\eqref{def:xj1},\eqref{def:generalGamma}}\ge  \beta   \Gamma_{j-1}(z_j)+\frac{1}{2\lambda}\|z_j-z_0\|^2 \nn \\
     & \overset{\beta <1}\ge \beta  \left[\Gamma_{j-1}(z_j)+ \frac{1}{2\lambda}\|z_j-z_0\|^2 \right]
     \stackrel{\eqref{def:xj1}}= \beta \Gamma_{j-1}^\lam(z_j) \nn \\
    &\ge  \beta \left[\Gamma_{j-1}^\lam(z_{j-1}) + \frac{1}{2\lam}\|z_j-z_{j-1}\|^2\right] 
 \label{eq:1},
\end{align}
where the last inequality follows from \eqref{def:xj1} with $j$ replaced by $j-1$ and the fact that $ \Gamma_{j-1}^\lam$ is $\lam^{-1}$-strongly convex.
Rearranging \eqref{eq:1}, taking the expectation of the resulting inequality, we have
\begin{align*}
\Exp{\Gamma_j^{\lam}(z_j)- \beta \Gamma_{j-1}^\lam(z_{j-1})} &\overset{\eqref{eq:1}}{\ge} (1-\beta)\Exp{\ell(z_j,z_{j-1};\xi_{j-1})} + \frac{\beta}{2\lam}\Exp{ \|z_j-z_{j-1}\|^2}\\
&\overset{\eqref{eq:basic}}{\ge}  (1-\beta)\Exp{\phi(z_j)} - 2(1-\beta)M \sqrt{\Exp{\|z_j-z_{j-1}\|^2}} +\frac{\beta}{2\lam} \Exp{\|z_j-z_{j-1}\|^2}
\end{align*}
where the second inequality is due to \eqref{eq:basic}.
Minimizing the right-hand side in the above inequality with respect to $\sqrt{\Exp{\|z_j-z_{j-1}\|^2}}$, we obtain
\begin{equation} \label{ineq:inter}
    \Exp{\Gamma_j^{\lam}(z_j) } \ge \beta
\Exp{\Gamma_{j-1}^\lam(z_{j-1})} 
+ \Exp{(1-\beta)\phi(z_j)} - \frac{2 \lam(1-\beta)^2 M^2}{\beta}. 
\end{equation}
Using the definitions of $\alpha_j$ and $u_j$ in \eqref{def:delta} and \eqref{def:u}, respectively, \eqref{eq:unbiasedzj}, and the above inequality, we conclude that
\begin{align*}
\alpha_j +   \frac{2\lam M^2(1-\beta)}{\beta} &\overset{\eqref{def:delta}} = \Exp{u_j - \Gamma_j^{\lam}(z_j) }  
\stackrel{\eqref{def:u},\eqref{eq:unbiasedzj}}\le\Exp{\beta u_{j-1}+(1-\beta)\phi(z_j)  -\Gamma_j^{\lam}(z_j)}  \\
& \overset{\eqref{ineq:inter}
}\le\beta  \Exp{ u_{j-1}- \Gamma_{j-1}^{\lam}(z_{j-1})} + \frac{2 \lam(1-\beta)^2 M^2}{\beta}\\
&\overset{\eqref{def:delta}} = \beta \alpha_{j-1} + 2\lam M^2(1-\beta)+ \frac{2 \lam(1-\beta)^2 M^2}{\beta},
\end{align*}
and hence that \eqref{ineq:alphaj} holds.
\end{proof}

\vgap

We are ready to present the proof of Proposition \ref{lem:t1-2}.

\vgap

\noindent
{\bf Proof of Proposition \ref{lem:t1-2}:}
It follows from \eqref{def:xj1} and the fact that the objective function of \eqref{def:xj1} is $\lam^{-1}$-strongly convex that 
	\begin{equation}\label{eq:2}
	\Gamma_I^\lam(z_I) + \frac{1}{2 \lam}\|x_*- z_I\|^2 
	\le  \Gamma_I^\lam(x_*) \stackrel{\eqref{def:xj1}}= \Gamma_I(x_*) + \frac{1}{2 \lam}\|x_*- z_0\|^2.
	\end{equation}
The above inequality, the definition of $\alpha_I$ in \eqref{def:delta},  and Lemma \ref{lem:102}, then imply  that
	\begin{align}
	      \Exp{u_I} &- \Exp{\Gamma_I(x_*)}  - \frac1{2\lam}\Exp{\|z_0 -x_*\|^2}+ \frac{1 }{2 \lam}\Exp{\| z_I - x_*\|^2} 
       \overset{\eqref{eq:2}}\le \Exp{u_I-\Gamma_I^\lam(z_I)} \nn \\
        &\stackrel{\eqref{def:delta}} = \alpha_I +
            \frac{2\lam M^2(1-\beta)}{\beta} \overset{\eqref{ineq:alphaj}}\le \alpha_1 \beta^{I-1}   +
      \frac{2\lam M^2(1-\beta)}{\beta}
	     \overset{\eqref{ineq:alpha1}}{\le} 2\lam M^2\left(\beta^{I-1}  + \frac{1-\beta}{\beta} \right). \label{ineq:core}
	\end{align}
 Note that  Lemma \ref{lem:tauI} with $C =I+1$ implies that 
\begin{equation}\label{obs1}
\beta^{I-1} \stackrel{\eqref{ineq:beta-1}}\le 9 \beta^{I+1}   \stackrel{\eqref{ineq:beta-1}}\le \frac{9}{I+1} \le \frac{18\log(I+1)}{I+1},
\end{equation}
where the last inequality is due to $\log (I+1) \ge 1/2$ for every $I \ge 1$, and 
\begin{equation}\label{obs2}
\frac{1-\beta}{\beta} \stackrel{\eqref{def:beta}}= \frac{2\log (I+1)}{I+1-\log (I+1)} \le \frac{2 \log (I+1)}{I+1-(I+1)/2} = \frac{4\log (I+1)}{I+1},
\end{equation}
where the last inequality is due to $\log (I+1) \le (I+1)/2$ for every $I \ge 1$.
    Plugging observations \eqref{obs1} and \eqref{obs2} into \eqref{ineq:core}, we conclude that
\[
 \Exp{u_I} - \phi_*  - \frac1{2\lam}\Exp{\|z_0 -x_*\|^2}+ \frac{1 }{2 \lam}\Exp{\| z_I - x_*\|^2} \stackrel{\eqref{ineq:core}}\le  \frac{44\lam M^2 \log (I+1)}{I+1} + \Exp{\Gamma_I(x_*)} - \phi_*.
	\]
Therefore, \eqref{keyrecu} immediately follows from the definition of ${\cal N}(x_*;\Gamma_I)$ in \eqref{def:Nj}.
\QEDA

\section{Numerical Experiments}\label{sec:numerics}

This section
benchmarks the numerical results
of two variants of the S-Max1C method
of Section~\ref{sec:method}
and two variants of the M-Max1C method
of Section \ref{sec:multistage}
against RSA from \cite{nemjudlannem09}, DA from \cite{nesterov2009primal},  and SCPB  from \cite{liang2024single}.
The two variants of each method differ only in the choice of the index set $B$
used in the S-Max1C subroutine, namely $B=\{1\}$ for the first variant, or $B$ defined
in~\eqref{B:2_i} for the second variant.

These comparisons are made on one stochastic programming problem studied in the numerical experiments of \cite{liang2024single} and four real-world applications.
The implementations use MOSEK~10.2%
\footnote{\url{https://docs.mosek.com/latest/toolbox/index.html}}
to generate the stochastic oracles $s(x,\xi)$ and to solve the subproblem~\eqref{def:xj1}.
The MATLAB and Python implementations of the numerical experiments are publicly available on GitHub%
\footnote{\url{https://github.com/Honghao-Zhang1/pace-stochastic-multicut}}.
The results reported in this paper were obtained using the MATLAB implementation, while the Python implementation is provided as an open-source alternative to facilitate reproducibility. The computations are performed on PACE \footnote{https://pace.gatech.edu/} with Dual Intel Xeon Gold 6226 CPUs @ 2.7 GHz (24 cores/node).

Now we start by describing the methods. First, note that for each problem, $M$ was estimated
as for RSA, i.e.,  taking the maximum of $\|s(\cdot,\cdot)\|$ over 10,000 calls to the stochastic oracle at randomly generated feasible solutions. Second, all the methods are run for 200 and 1000 iterations. 


\par {\textbf{RSA:}}
Given an iterate $x_t$ and a stochastic subgradient $g_t$, RSA described in Section~2.2 of~\cite{nemjudlannem09}
updates according to
\[
x_{t+1} = \argmin_{x \in X} \left\{ \IInner{g_t}{x}+\frac{1}{2\gamma_t}\|x-x_t\|^2\right\}.
\]
The output is 
$\displaystyle  x^N=\frac{1}{N} \sum_{i=1}^N x_i$
where $x_i$
is computed at the $i$-th iteration using the constant
steps given in (2.23) of \cite{nemjudlannem09}
by 
\begin{equation}\label{def:rsagamma}
\gamma_t=\frac{C D}{M \sqrt{N}}
\end{equation}
where $D$ is the diameter of the feasible set $X$.
As in \cite{nemjudlannem09}, we take $C=0.1$
which was calibrated in \cite{nemjudlannem09}
using an instance of the stochastic 
utility problem. $N \in \{200,1000\}$.

\vspace*{0.2cm}
\par {\textbf{SCPB:}}
SCPB is as described in Section 6 of \cite{liang2024single} (denoted as SCPB1 in \cite{liang2024single})  and uses
parameters
$\theta$, $\tau$, 
$R$, and $\lambda$  given 
by 
\[
\theta = \frac{9}{K},
\;\tau=\frac{\theta K}{\theta K +1},\;
R=\frac{D}{M},\;
\lambda= 10 \frac{3 D}{M \sqrt{K}}.
\]
For each targeted  iteration count $N \in \{200,1000\}$, $K$ is set to be $N/10$.

\vspace*{0.2cm}
\par {\textbf{S-1C:}} 
S-1C is S-Max1C with $B = \{1\}$ and uses stepsize $\lam$ given by 
$
\lam = (10 \sqrt{I} D)/M$ where $I \in \{200,1000\}$.

\vspace*{0.2cm}
\par {\textbf{S-Max1C:}} 
S-Max1C with $B$ as in \eqref{B:2_i} uses stepsize $\lam$ given by 
$\lam = (10 \sqrt{I} D)/M $ where $I \in \{200,1000\}$.
\vspace*{0.2cm}
\par {\textbf{M-1C:}} 
M-1C is M-Max1C with $B=\{1\}$,  uses two phases $L=2$, and sets stepsize
$
\lam = (10 \sqrt{I} D)/(\sqrt{2}M)$  where  $I=N/2$ for total number of iterations $N \in \{200,1000\}$.

\vspace*{0.2cm}
\par {\textbf{M-Max1C:}} 
M-Max1C uses stepsize $\lam$ given by 
$\lam = (10 \sqrt{I} D)/(\sqrt{2}M)$ where  $I=N/2$ for total number of iterations $N \in \{200,1000\}$,
and sets $B = \{2^i:i \ge 0, 2^i \le \lfloor I/2\rfloor\}$ and $L=2$. 

\vspace*{0.2cm}
\par {\textbf{DA:}} DA is as described in \cite{nesterov2009primal} and updates as follows:
\[
 x_{k+1} = \argmin_{x \in X} \left\{ \IInner {\sum_{i=0}^k g_i}{x} +\frac{ \gamma_k}{2}\|x-x_0\|^2\right\}
\]
where $g_i $ is a stochastic subgradient of $f$ at $x_i$.
We choose the stepsize given by \cite{nesterov2009primal}  as
\begin{equation}\label{def:DAgamma}
\gamma_k = \frac{M \alpha_k} {C\sqrt{D}},
\end{equation}
where $C=10$, $\alpha_1=\alpha_0 =1$, and for all $k \ge 2$
\[
\alpha_k = \alpha_{k-1} + \frac{1}{\alpha_{k-1}}.
\]



\par \textbf{Notation in the tables.} In the following, we use the following notation:
\begin{itemize}
\item $n$ represents the design dimension of an instance;
\item $N$ denotes the sample size used to run the methods; this also corresponds to the number of iterations of RSA; 
\item {\tt{Obj}} refers to the empirical mean 
\begin{equation}\label{defempmean}
{\hat F}_T(x) := \frac{1}{T} \sum_{i=1}^T F(x, \xi_i)
\end{equation}
of $F$ at $x$ based on a sample 
$\xi_1, \ldots, \xi_T$ of $\xi$ with size $T$, which provides an estimate of $f(x)$. The 
empirical means are computed with $x$ being the final iterate output
by the algorithm and $T = 10^4$;
\item {\tt{CPU}} refers to the rounded CPU time in seconds;
\item {\tt{Std}} refers to the standard deviation of the observed objective values.
\end{itemize}

\subsection{Two-stage stochastic program}\label{subsec:sto2}
We consider a nonlinear two-stage stochastic program 
\begin{equation}\label{smdmodel21}
\left\{
\begin{array}{l}
\min \; c^\top x_1 + \mathbb{E}[\mathfrak{Q}(x_1,\xi)]\\
x_1 \in \mathbb{R}^n : \|x_1\|_2 \leq D 
\end{array}
\right. 
\end{equation}
where the cost-to-go function $\mathfrak{Q}(x_1,\xi)$ has nonlinear objective and constraint coupling functions and is given by
\begin{equation}\label{smdmodel22}
\mathfrak{Q}(x_1,\xi)=\left\{ 
\begin{array}{l}
\displaystyle  \min_{x_2 \in \mathbb{R}^n} \;\frac{1}{2}\left( \begin{array}{c}x_1\\x_2\end{array} \right)^\top \Big( \xi \xi^\top + \gamma_0 I_{2 n} \Big) \left( \begin{array}{c}x_1\\x_2\end{array} \right) + \xi^T \left( \begin{array}{c}x_1\\x_2\end{array} \right)  \\
\text{s.t.} \ \ \ \|x_2\|_2^2 + \|x_1\|_2^2 - R^2 \leq 0.
\end{array}
\right.
\end{equation}
This problem is an instance of SCCO \eqref{eq:ProbIntro}-\eqref{pbint2}, with $h$ being the indicator function of a compact convex set
$X$ with diameter
$D$.

Throughout this subsection, $\xi$ is generated as a Gaussian random vector in $\mathbb{R}^{2n}$  with means and standard deviations 
uniformly sampled at random  in 
$[-\chi,\chi]$  and $[0,\chi]$, respectively, for some $\chi >0$ specified in corresponding problem instances. The components of $c$ are generated uniformly at random from $[-1,1]$.  Parameter $\gamma_0 $ is set to $2$. 

\begin{table}[H]
\scriptsize
\centering
\setlength{\tabcolsep}{3pt}
\begin{tabular}{|c|c|ccc|ccc|ccc|ccc|}
\hline
\multicolumn{2}{|c|}{} &
\multicolumn{6}{|c|}{$(D,R,\chi)=(2,4,5)$} &
\multicolumn{6}{|c|}{$(D,R,\chi)=(50,100,2)$} \\
\hline
\multicolumn{2}{|c|}{} &
\multicolumn{3}{|c|}{$C_1: n=100$} &
\multicolumn{3}{|c|}{$C_2: n=200$} &
\multicolumn{3}{|c|}{$C_3: n=100$} &
\multicolumn{3}{|c|}{$C_4: n=200$} \\
\hline
ALG. & N& {\tt{Obj}}& {\tt{Std}} &{\tt{CPU}} &
{\tt{Obj}}& {\tt{Std}} &{\tt{CPU}} &
{\tt{Obj}}& {\tt{Std}} &{\tt{CPU}} &
{\tt{Obj}}& {\tt{Std}} &{\tt{CPU}}\\
\hline
RSA&200 & -5.657 & 0.57 & 7.8 & -9.406 & 0.52 & 42.8 & -7.204 & 0.79 & 8.1 & -14.581 & 1.28 & 43.4 \\
&1000& -7.116 & 0.57 & 39.4 & -11.531 & 0.51 & 213.0 & \textbf{-8.275} & 0.84 & 40.4 & -16.502 & 1.32 & 216.6 \\
\hline
DA&200 & \textbf{-8.505} & 0.97 & 7.8 & -0.321 & 0.44 & 40.0 & 8.045 & 2.08 & 8.1 & 7.500 & 2.08 & 43.3 \\
&1000&\textbf{-8.538} & 0.87 & 38.8 & -0.638 & 0.20 & 200.6 & 1.625 & 1.57 & 40.3 & -6.112 & 1.75 & 216.2 \\
\hline
S-1C&200 & -7.838 & 0.62 & 7.7 & \textbf{-12.773} & 0.57 & 42.6 & 0.878 & 1.36 & 8.1 & -6.532 & 1.23 & 43.7 \\
&1000& -7.838 & 0.62 & 38.4 & \textbf{-12.774} & 0.57 & 212.5 & -6.110 & 1.00 & 40.3 & -14.171 & 1.46 & 218.1 \\
\hline
S-Max1C&200 & -7.837 & 0.62 & 8.0 & \textbf{-12.773} & 0.57 & 43.3 & \textbf{-7.436} & 0.94 & 8.4 & \textbf{-15.709} & 1.26 & 44.3 \\
&1000& -7.838 & 0.62 & 40.3 & \textbf{-12.774} & 0.57 & 218.1 & -8.201 & 0.83 & 42.3 & \textbf{-16.586} & 1.32 & 220.8 \\
\hline
SCPB&200 & -7.837 & 0.62 & 8.3 & -12.772 & 0.57 & 46.4 & -3.137 & 1.61 & 8.8 & -10.276 & 2.69 & 47.8 \\
&1000& -7.838 & 0.62 & 41.1 & -12.773 & 0.57 & 226.4 & -7.121 & 0.96 & 43.2 & -15.273 & 1.33 & 232.4 \\
\hline
M-1C&200 & -7.838 & 0.62 & 8.1 & \textbf{-12.773} & 0.57 & 44.7 & 0.878 & 1.36 & 8.5 & -6.532 & 1.23 & 45.9 \\
&1000& -7.838 & 0.62 & 40.3 & \textbf{-12.774} & 0.57 & 223.1 & -6.110 & 1.00 & 42.3 & -14.171 & 1.46 & 229.0 \\
\hline
M-Max1C&200 & -7.836 & 0.63 & 8.4 & \textbf{-12.773} & 0.58 & 45.5 & 7.876 & 8.20 & 8.8 & 0.173 & 10.99 & 46.5 \\
&1000& -7.838 & 0.62 & 42.3 & \textbf{-12.774} & 0.57 & 228.0 & -4.812 & 2.30 & 44.4 & -14.064 & 3.12 & 231.8 \\
\hline
\end{tabular}
\caption{Average performance over all seeds on the two-stage stochastic program \eqref{smdmodel21}-\eqref{smdmodel22}. Best (lowest) average objective values in each column are highlighted in bold.}
\label{tabsecondpb:cpuobj3}
\end{table}

We run RSA, SCPB, DA, and two versions of S-Max1C and M-Max1C on 
instances $C_1$-$C_4$. For  instances $C_1$ and $C_2$ (resp., $C_3$ and $C_4$), $(D,R,\chi)$ is set to $(2,4,5)$ (resp., $(50,100,2)$).
Each method is executed over $30$ independent runs, and the results are reported in Table \ref{tabsecondpb:cpuobj3}.
Bold numbers highlight the best algorithm in terms of objective value for each target iteration ($200$ and $1000$).

We now make several comments about Table \ref{tabsecondpb:cpuobj3}.  First, S-1C and S-Max1C achieve the best performance in most cases. Second, S-Max1C is comparable to S-1C on instances $C_1$ and $C_2$, and substantially outperforms S-1C on $C_3$ and $C_4$.
This performance demonstrates the superiority of the Max1C model over the one-cut model (S1) (see \eqref{eq:Gammaj}) on these instances. However, in contrast, M‑Max1C performs poorly on instances C3 and C4 compared with the M-1C.

\subsection{Real-world problems}

This subsection presents numerical experiments on four real-world applications, including motor freight carrier operations \cite{mak1999monte}, aircraft allocation \cite{Dantzig1963}, telecommunications network design \cite{sen1994network}, and cargo flight scheduling  \cite{mulvey1995new}, which are all modeled as a two-stage stochastic linear program with recourse.  In the standard form, it can be formulated as
\[
\min_{x} \; c^{\top}x + \mathcal{Q}(x)
\quad \text{subject to} \quad
Ax = b,\; x \ge 0,
\]
where \(\mathcal{Q}(x) := \mathbb{E}[Q(x,\xi)]\) denotes the expected recourse function, \(Q(x,\xi)\) is defined as the optimal value of the second-stage problem
\[
Q(x,\xi) := \min_{y} \; q^{\top}y
\quad \text{subject to} \quad
Tx + Wy = r,\; y \ge 0,
\]
and $\xi = (q,r,T,W)$ collects the stochastic data of the second-stage problem.
The input datasets used in our experiments, provided in SMPS format, are \textsc{20term}, \textsc{gbd}, \textsc{ssn}, and \textsc{storm}, available at \url{https://pages.cs.wisc.edu/~swright/stochastic/sampling/}.

Since this paper focuses on single-phase algorithms, we benchmark four single-phase SA methods, namely RSA, DA, S-1C, and S-Max1C, across the four real-world applications. To provide a comprehensive comparison, we evaluate each method under four different stepsize settings and report the best-performing result in terms of the average objective value. Specifically, for both S-1C and S-Max1C, we consider four candidate stepsizes of the form
\[
\lambda = \frac{C\sqrt{I}\,D}{M}, \quad C \in \{0.0001,\; 0.01,\; 1,\; 10\},
\]
and report the best result obtained for each instance. Similarly, for RSA (with parameter $C$ in \eqref{def:rsagamma}) and DA (with parameter $C$ in \eqref{def:DAgamma}), we test four values $
\{0.1,\; 1,\; 5,\; 10\}$,
and present the best outcome among them.



Tables  \ref{tabsecondpb:cpuobj20term}, \ref{tabsecondpb:cpuobjgbd}, \ref{tabsecondpb:cpuobjSSN}, and \ref{tabsecondpb:cpuobjstorm} below report numerical results comparing the performance of the four methods on the four datasets  \textsc{20term},  \textsc{gbd},  \textsc{ssn}, and \textsc{storm}, respectively. Each method is executed over $30$ independent runs for two different sample sizes $200$ and $1000$.

The \textsc{20term} problem \cite{mak1999monte} models a motor freight carrier's operations with stochastic shipment demands. The objective is to position a fleet and route vehicles through a network to satisfy point-to-point demands while penalizing unmet demand and enforcing end-of-day fleet balance. The average performance results are reported in Table~\ref{tabsecondpb:cpuobj20term}.

\begin{table}[htb]
\footnotesize
\centering
\begin{tabular}{|c|c|ccc|}
\hline
ALG. & $N$ & {\tt{Obj}} & {\tt{Std}} & {\tt{CPU}} \\
\hline
\multirow[t]{2}{*}{RSA}
&200  & 269620 & 541.71 & 3.8  \\
&1000 & 259650 & 290.92 & 18.8 \\
\hline
\multirow[t]{2}{*}{DA}
&200  & 254500 & 280.52 & 4.1  \\
&1000 & 254430 & 282.14 & 20.8 \\
\hline
\multirow[t]{2}{*}{S-1C}
&200  & 254500 & 277.31 & 3.9  \\
&1000 & 254460 & 280.70 & 19.3 \\
\hline
\multirow[t]{2}{*}{S-Max1C}
&200  & 254500 & 278.91 & 4.0 \\
&1000 & 254460 & 280.12 & 20.2 \\
\hline
\end{tabular}
\caption{Performance of RSA, DA, S-1C, and S-Max1C on the \textsc{20term} problem with 30 runs.}
\label{tabsecondpb:cpuobj20term}
\end{table}

The \textsc{gbd} problem \cite{Dantzig1963} is an aircraft allocation model that assigns multiple aircraft types to routes to maximize expected profit under uncertain passenger demand. The formulation accounts for operating costs and penalties for unmet demand, with uncertainty represented through a large set of demand scenarios. The average performance results are reported in Table~\ref{tabsecondpb:cpuobjgbd}.

\begin{table}[H]
\footnotesize
\centering
\begin{tabular}{|c|c|ccc|}
\hline
ALG. & $N$ & {\tt{Obj}} & {\tt{Std}} & {\tt{CPU}}\\ 
\hline
\multirow[t]{2}{*}{RSA} 
&200  & 1668.9 & 20.58 & 1.7 \\
&1000 & 1665.6 & 20.31 & 8.4 \\
\hline
\multirow[t]{2}{*}{DA} 
&200  & 1661.8 & 20.30 & 1.7 \\
&1000 & 1661.7 & 20.25 & 8.3 \\
\hline
\multirow[t]{2}{*}{S-1C} 
&200  & 1661.5 & 20.25 & 1.7 \\
&1000 & 1661.5 & 20.25 & 8.5 \\
\hline
\multirow[t]{2}{*}{S-Max1C} 
&200  & 1661.5 & 20.25 & 1.8 \\
&1000 & 1661.5 & 20.25 & 9.1 \\
\hline
\end{tabular}
\caption{Performance of RSA, DA, S-1C, and S-Max1C on the \textsc{gbd} problem with 30 runs.}
\label{tabsecondpb:cpuobjgbd}
\end{table}


The \textsc{ssn} problem \cite{sen1994network} models a telecommunications network design task in which random service requests must be routed with sufficient capacity. The goal is to determine capacity expansions that minimize the expected rate of unmet demand under stochastic scenarios. The average performance results are reported in Table~\ref{tabsecondpb:cpuobjSSN}.
\begin{table}[H]
\footnotesize
\centering
\begin{tabular}{|c|c|ccc|}
\hline
ALG. & $N$ & {\tt{Obj}} & {\tt{Std}} & {\tt{CPU}}\\ 
\hline
\multirow[t]{2}{*}{RSA} 
&200  & 10.2578 & 0.54 & 2.9 \\
&1000 & 10.0373 & 0.52 & 14.7 \\
\hline
\multirow[t]{2}{*}{DA} 
&200  & 9.8456 & 0.53 & 3.1 \\
&1000 & 9.8372 & 0.53 & 15.7 \\
\hline
\multirow[t]{2}{*}{S-1C} 
&200  & 9.8370 & 0.53 & 3.3 \\
&1000 & 9.8364 & 0.53 & 16.0 \\
\hline
\multirow[t]{2}{*}{S-Max1C} 
&200  & 9.8364 & 0.52 & 3.2 \\
&1000 & 9.8364 & 0.52 & 16.3 \\
\hline
\end{tabular}
\caption{Performance of RSA, DA, S-1C, and S-Max1C on the \textsc{ssn} problem with 30 runs.}
\label{tabsecondpb:cpuobjSSN}
\end{table}


The \textsc{storm}  problem \cite{mulvey1995new} is based on a cargo flight scheduling application. The aim is to plan cargo-carrying flights
over a set of routes in a network, where the amounts of cargo are uncertain. The average performance results are reported in Table~\ref{tabsecondpb:cpuobjstorm}.
\begin{table}[H]
\footnotesize
\centering
\begin{tabular}{|c|c|ccc|}
\hline
ALG. & $N$ & {\tt{Obj}} & {\tt{Std}} & {\tt{CPU}}\\ 
\hline
\multirow[t]{2}{*}{RSA} 
&200  & 1.5500$\times 10^7$ & 8521.3 & 5.34 \\
&1000 & 1.5500$\times 10^7$ & 8523.1 & 26.10 \\
\hline
\multirow[t]{2}{*}{DA} 
&200  & 1.5500$\times 10^7$ & 8534.0 & 5.86 \\
&1000 & 1.5500$\times 10^7$ & 8533.4 & 28.96 \\
\hline
\multirow[t]{2}{*}{S-1C} 
&200  & 1.5500$\times 10^7$ & 8534.5 & 4.17 \\
&1000 & 1.5500$\times 10^7$ & 8535.8 & 21.07 \\
\hline
\multirow[t]{2}{*}{S-Max1C} 
&200  & 1.5500$\times 10^7$ & 8516.9 & 5.44 \\
&1000 & 1.5500$\times 10^7$ & 8545.0 & 32.51 \\
\hline
\end{tabular}
\caption{Performance of RSA, DA, S-1C, and S-Max1C on the \textsc{storm} problem with 30 runs.}
\label{tabsecondpb:cpuobjstorm}
\end{table}

We finally conclude this subsection with several remarks. First, S-1C and S-Max1C generally outperform RSA and achieve performance comparable to DA.
Second, although S-Max1C incurs a slightly higher per-iteration computational cost compared to the other algorithms due to its step~2, S-Max1C and S-1C generally exhibit comparable standard
deviations.

\section{Concluding Remarks}\label{sec:conclusion}

This paper studies multi-cut SA methods for solving SCCO \eqref{eq:ProbIntro}
that relax the  condition \eqref{ineq:common}. It proposes the generic S-CP framework and, specifically studies one of its instances, namely, the S-Max1C method, which
is based on a cutting-plane model \eqref{def:Gammaj} lying  between the multi-cut model \eqref{multicut} and the one-cut model \eqref{onecut}. It is shown that S-Max1C has the convergence rate $\tilde {\cal O}(1/\sqrt{I})$, which is the same as other one-cut SA methods up to a logarithmic term. Leveraging a warm-start approach,  a multi-phase version of S-Max1C, i.e., the M-Max1C method, is developed and shown to have the same convergence rate (in terms of $I$) as S-Max1C. Computational results demonstrate that both S-Max1C and M-Max1C are generally comparable to  and sometimes outperform standard SA methods  in all instances considered.

We provide some remarks and possible extensions of this paper. First, the assumption that the function $F(x,\xi)$ is convex in $x$ (see (A2)) can actually be removed at the expense of a technically more involved proof of Lemma \ref{lem:basicLQ}, but we make this assumption for the sake of simplicity.
Second, the convergence rate of M-Max1C in Theorem~\ref{cor:cmplx1-practical2} is not optimal in terms of both $I$ and $L$, and hence it would be interesting to develop a multi-phase multi-cut SA method that enjoys the optimal rate ${\cal O} (1/\sqrt{LI})$.
Third, the noise for the max-one-cut model, ${\cal N}(\cdot;\Gamma_j) = {\cal O}(1/\sqrt{j})$, may be attributed to the intrinsic nonsmoothness of the pointwise maximum function. It is worth investigating whether using smooth approximations of the max function can help improve the order of the noise behavior. One potential approach could be to replace the maximum function in the multi-cut model with the LogSumExp function. Finally, it is of practical interest to explore more possibilities of selecting Max1C index set $B$ beyond the two options used in S-1C and S-Max1C.

\bibliographystyle{plain}
 \bibliography{ref}
 \appendix

\section{Technical Results}\label{app:technical}

\begin{lemma}\label{lem:ind}
    Assume $\{Z_i\}_{i=1}^n$ are independent random variables such that $\Exp{Z_i}=0$ for every $i=1,\ldots,n$. Then,
\[
\Exp{ \left (\sum_{i=1}^n Z_i \right)^2 }
= \sum_{i=1}^n \Exp{Z_i^2}.
\]
\end{lemma}

\begin{lemma}\label{lem:variance}
    For some finite index set $B$ and scalar $ \sigma_X \ge 0$, assume that $\{Y_k\}_{k \in B}$ and $\{X_k\}_{k \in B}$ are families of real-valued random variables such that
    \begin{equation}\label{eq:condition}
        Y_k \le X_k , \quad \Exp{X_k} = 0, \quad \Var(X_k) \le \sigma_X^2,
    \quad \forall k \in B.
    \end{equation}
    Then, we have 
    \begin{equation}\label{ineq:E-max}
        \Exp{\max_{k \in B} Y_k} \le 2\sigma_X \sqrt{|B|-1}.
    \end{equation}
\end{lemma}

\begin{proof}
First, observe that the assumption
implies that $\Exp{X_k^2} \le \sigma_X^2$ for every $k \in B$.
Set $\bar X=X_i$ for some fixed $i \in B$.
The above observation then implies that 
\begin{equation}\label{eq:vr}
    \Exp{(X_k-\bar X)^2 } \le 2\Exp{ X_k^2 } + 2\Exp{\bar X^2 } \le 4\sigma_X^2.
\end{equation}
where  the first inequality is due to $(a_1+a_2)^2\le 2a_1^2+2a_2^2$ for every $a_1,a_2\in \R$.
It thus follows from the fact that $\|\cdot\|_\infty \le \|\cdot\|_2$ that
    \begin{align*}
        \Exp{\max_{k \in B} Y_k} 
 &\le\Exp{\max_{k \in B} X_k} = \Exp{\left(\max_{k \in B}   X_k\right) - X_i  } 
 = \Exp{ \max_{k \in B} \{  X_k - \bar X \} } \le \Exp{ \left( \sum_{ k \in B} (X_k - \bar X)^2 \right)^{1/2} }
 \\
&= \Exp{ \left( \sum_{ k \in B\setminus \{i\} } (X_k - \bar X)^2 \right)^{1/2} }
\le \left( \Exp{ \sum_{ k \in B\setminus \{i\} } (X_k - \bar X)^2  }\right)^{1/2}
\overset{\eqref{eq:vr}}\le \sqrt{4 \sigma_X^2 \left(|B|-1 \right) },
    \end{align*}
  and hence \eqref{ineq:E-max} holds.
\end{proof}

\begin{lemma}\label{lem:tauI}
Let  $C \ge 2$ be given and define $\beta := (C-\log C)/(C+\log C)$.
Then,
\begin{equation}\label{ineq:beta-1}
     \beta^C \le \frac{1}{C}, \quad 
\beta \ge \frac{1}{3}.
\end{equation}
\end{lemma}

\begin{proof}
 We first prove the first inequality. Indeed, using the definition of $\beta$ and the fact that $\log x \le x -1 $ for any $x>0$, we have
\[
    \beta^C = e^{C \log \beta} \le e^{C(\beta -1)} = e^{-\frac{2 C \log C}{C + \log C}} \le e^{- \log C} =\frac{1}{C},
\]
where the last inequality follows from the fact that $\log C \le C$.
The second inequality follows from the fact that
the assumption that $C \ge 2$ implies that $C/\log C \ge 2$ and the fact that
the function
$t \in [1,\infty) \mapsto (t-1)/(t+1)$ is increasing.
\end{proof}


\begin{lemma}\label{lem:Nj}
    For every $j\ge 1$, consider ${\cal N} (\cdot;\Gamma_j)$ in \eqref{def:Nj} computed for $\Gamma_j$ as in \eqref{eq:multicut}, then for every $u \in \dom h$, we have
\begin{equation}\label{ineq:badbd1}
{\cal N}(u;\Gamma_j) \le 2\sigma(u)\sqrt{j-1} ,
\end{equation}
where $\sigma(u)$ is defined as in \eqref{def:sigma_u}.
\end{lemma}

\begin{proof}
In view of the definition of $\sigma(u)$ given in \eqref{def:sigma_u}, relation \eqref{ineq:linear} implies that for given $u \in \dom h$, the random variables $\{X_k\}_{k \in B}$ and $\{Y_k\}_{k \in B}$,  the scalar $\sigma_X$, and the index set $B$, defined as
\[
X_k = \Phi(u,\xi_k) - \phi(u), \quad
Y_k =\ell(u,z_k;\xi_k) - \phi(u),\quad \sigma_X = \sigma(u), \quad B=\{0,1,\ldots,j-1\},
\]
satisfy \eqref{eq:condition}.
Hence, it follows from the conclusion of 
 Lemma \ref{lem:variance} and  the definition of $\Gamma_j$ in \eqref{eq:multicut} that 
\[
\Exp{\Gamma_j(u)}-\phi(u) \stackrel{\eqref{eq:multicut}}= \Exp{\max_{0\le k \le j-1} \left\{\ell(u,z_k;\xi_k)-\phi(u)\right\}} \stackrel{\eqref{ineq:E-max}}\le 2\sigma(u)\sqrt{j-1}.
\]
Finally, inequality \eqref{ineq:badbd1} follows from the above inequality and the definition of ${\cal N}(u;\Gamma_j)$ in \eqref{def:Nj}. 
\end{proof}







\end{document}

%% file: def.tex
\usepackage{amsmath}
\usepackage{amsfonts}
\usepackage{latexsym}
\usepackage{amssymb}

\newtheorem{thm}{Theorem}[section]
\newtheorem{theorem}[thm]{Theorem}

\newtheorem{lemma}[thm]{Lemma}

\newtheorem{proposition}[thm]{Proposition}

\newtheorem{remark}[thm]{Remark}

\newcommand{\Exp}[1]{
\mathbb{E}\left[#1 \right]
}
\newcommand{\beq}{\begin{equation}}
\newcommand{\eeq}{\end{equation}}
\newcommand{\beqa}{\begin{eqnarray}}
\newcommand{\eeqa}{\end{eqnarray}}
\newcommand{\beqas}{\begin{eqnarray*}}
\newcommand{\eeqas}{\end{eqnarray*}}
\newcommand{\bi}{\begin{itemize}}
\newcommand{\ei}{\end{itemize}}

\newcommand{\vgap}{\vspace{.1in}}

\newcommand{\nn}{\nonumber}

\newcommand{\R}{\mathbb{R}}

\newcommand{\lam}{{\lambda}}

\newcommand{\inner}[2]{\langle #1,#2\rangle}

\newcommand{\IInner}[2]{\left \langle #1\,,#2 \right \rangle}

\newcommand{\argmin}{\mathrm{argmin}\,}

\newcommand{\dom}{\mathrm{dom}\,}

\newcommand{\bConv}[1]{\overline{\mbox{\rm Conv}}\,(\R^{#1})}